\documentclass[11pt,english]{article}
\usepackage[T1]{fontenc}

\usepackage{graphicx, subfigure}
\expandafter\def\csname ver@subfig.sty\endcsname{}
\usepackage[latin9]{inputenc}
\usepackage{geometry}
\usepackage{amsmath}
\usepackage{appendix}
\usepackage{amssymb}
\usepackage{esint}
\usepackage{makecell,multirow,diagbox}
\usepackage{mathtools}
\usepackage{epstopdf}
\usepackage[latin9]{inputenc}
\usepackage{pict2e}
\usepackage{xcolor}
\usepackage{colortbl,booktabs}
\usepackage{stmaryrd}
\usepackage{esint}
\usepackage[all,cmtip]{xy}
\usepackage{mathtools}
\usepackage{float}
\usepackage{yfonts}
\usepackage{setspace}
\usepackage{authblk}
\usepackage{amsthm}
\usepackage{mathrsfs}
\usepackage{stmaryrd}
\usepackage{bbold}
\usepackage{color}
\usepackage{tikz}
\usepackage[all,cmtip]{xy}
\usepackage{algorithm}  
\usepackage{algorithmic}

\theoremstyle{plain}
\theoremstyle{plain}
\newtheorem{remark}{Remark}
\usepackage{color}
\definecolor{marin}{rgb} {0., 0.3, 0.7}
\definecolor{rouge}{rgb} {0.8, 0., 0.}
\definecolor{sepia}{rgb} {0.8, 0.5, 0.}
\usepackage[colorlinks,citecolor=sepia,linkcolor=marin,bookmarksopen,bookmarksnumbered]{hyperref}

\theoremstyle{definition}

\DeclareSymbolFont{largesymbol}{OMX}{yhex}{m}{n}
\DeclareMathAccent{\Widehat}{\mathord}{largesymbol}{"62}

\makeatletter

\usepackage{babel}
\usepackage{subfig}
\graphicspath{{Dirac/}}
           
\makeatother
\DeclareMathSizes{14}{14}{9.8}{7}

\usepackage{lipsum}

\begin{document}

\title{\bf{Energy conserving particle-in-cell methods for relativistic Vlasov--Maxwell equations of laser-plasma interaction}}
\date{}
\author{\small \textbf{Yingzhe Li}}
\affil{Max Planck Institute for Plasma Physics, Boltzmannstrasse 2, 85748 Garching, Germany\\
yingzhe.li@ipp.mpg.de}
\maketitle
\begin{abstract}
Energy conserving particle-in-cell schemes are constructed for a class of reduced relativistic Vlasov--Maxwell equations of laser-plasma interaction. Discrete Poisson equation is also satisfied by the numerical solution. Specifically, distribution function is discretized using particle-in-cell method,  discretization of electromagnetic fields is done using compatible finite element method in the framework of finite element of exterior calculus, and time discretization used is based on discrete gradient method combined with Poisson splitting. Numerical experiments of parametric instability are done to validate the conservation properties and good long time behavior of the numerical methods constructed. 
\end{abstract}

\setcounter{tocdepth}{2} 


\section{Introduction}
Laser-plasma interaction is an important physical concept in the fields of inertial fusion confinement and plasma based electron accelerator schemes,  which include a lot of complex physical processes when strong lasers are injected into plasmas. When the plasma density is very high and particles are accelerated by the lasers to high speeds, the relativistic and quantum effects (such as spin effects) are unignorable. There are extensive theoretical, experimental, and numerical works about laser-plasma interaction. For example, in~\cite{Shengzm} the acceleration of electrons in plasma by two counter-propagating laser pulses is discussed, and numerical simulations are done for the interaction between spin-polarized electrons beams and strong laser pulses in~\cite{Wen}.  

Kinetic equations are adopted by the laser-plasma community for theoretical and numerical explanations.  As lasers usually propagate along fixed directions,  the models with lower dimensions reduced from three dimensional Vlasov--Maxwell equations can be used. In this sprit, there are one and two dimensional reduced laser-plasma models proposed in the literature~\cite{D1, D2}, in which the reduction relies on the conservation of the canonical momentum of particles. There are a lot of existing theoretical and numerical works about these laser-plasma models, such as~\cite{THEO1, THEO2, NUMA}, in which existence of mild and global solutions are done, also an error estimate result of a semi-Lagrangian method is given. To include spin effects, a set of kinetic equations is introduced recently and detailed in~\cite{spin1, spin2, spin3}.  And in~\cite{LPFEEC} a structure-preserving method for non-relativistic Vlasov--Maxwell equations with spin effects is introduced based on the geometric structures proposed in~\cite{LPS, LP}. In this work, we focus on the fully relativistic case. 

There are mainly two classes of methods for solving kinetic models in plasma physics, the grid based method and particle-in-cell method~\cite{pic1, pic2}. Grid based method includes for instances semi-Lagrangian method~\cite{semi1}, discontinuous Galerkin method~\cite{DG}, and so on.  When the dimension of phase-space and the domain scale of simulation are very large, grid based method are relatively costly, but without numerical noise, which decreases as $\mathcal{O}(1/\sqrt{N_p})$ about particle number $N_p$ in particle in cell method. The advantage of particle in cell method is the efficiency especially for high dimensional models. The reason we choose particle in cell method in this paper is that the spin variable in our kinetic model is sampled on the unit sphere, and thus is more suitable to discrete using particles. 

Our discretization follows the recent trend of structure-preserving methods~\cite{Feng, HLW}, which have been proposed with the purpose of preserving the intrinsic properties inherited by the given system and thus have long term stability and accuracy. In plasma physics, some structure-preserving methods~\cite{Qin2, Qin3, Qin4, Qin5, GEMPIC, curviSISC, Morrison, LPFEEC, martin} have been proposed for Vlasov type equations. In these works, space discretizations are done in the framework of finite element exterior calculus~\cite{FEEC} or discrete exterior calculus~\cite{DEC}, after which (time-continuous) finite dimensional Poisson systems (non-canonical Hamiltonian system) are derived. From~\cite{HLW}, we know that the only time discretization used to construct fully discrete structure-preserving methods for general non-canonical Hamiltonian systems is the so-called Hamiltonian splitting method~\cite{Nicolas, Qin4}, which requires each Hamiltonian subsystem obtained explicitly solvable, and thus do not work well for many complicated Hamiltonian models, especially when Hamiltonians are complicated. 

Under this circumstances, constructing methods preserving other theoretical properties, such as energy and constraints are meaningful for good long time simulations. As for the energy-conserving method, when the Hamiltonian is a quadratic function,  it can be conserved by the usual mid-point rule or Crank-Nicolson method. For more complicated Hamiltonians, discrete gradient method~\cite{DIS} has been proposed, which is adopted in this work. As mentioned in~\cite{2020en}, in which discrete gradient method is used to construct energy conserving schemes for non-relativistic Vlasov--Maxwell equations, when the Hamiltonian is a quadratic function, many existing discrete gradient methods will become mid-point rule method.  Another way to construct energy-conserving methods is the recently proposed so-called scalar auxiliary variable (SAV) approach~\cite{SAV1}, by which an equivalent new Hamiltonian could be conserved, while the original one is not conserved by the numerical scheme. As for relativistic Vlasov--Maxwell equations, a quadratic conservative finite difference method is proposed to conserved energy in~\cite{SOS}; an energy-conserving  finite difference method is proposed based on mid point rule and  Crank--Nicolson method in~\cite{Chacon}. An energy-conserving discontinuous Galerkin methods is proposed in~\cite{heyang}. An Eulerian conservative splitting scheme is proposed in~\cite{LP} based on Poisson structure of the system and has good long time behavior. The advantages of the numerical methods constructed in this work include: a) higher space accuracy can be obtained by increasing the degrees of basis functions of finite element spaces; b) there is no smoothing effect from discretizing particles, as delta functions are used rather than smoothed delta functions; c) energy is conserved and discrete Poisson equation is satisfied by the numerical solution as well; (d) the schemes can be extended to three dimensional case directly.

The paper is organized as follows. In section 2,  one and two dimensional laser-plasma models are introduced, specifically a Poisson bracket for the two dimensional case is proposed for the first time. In section 3, phase space discretizations are described, and finite dimensional Poisson systems with complicated (non-quadratic) Hamiltonians are derived. In section 4, energy conserving schemes are constructed using discrete gradient and Poisson splitting method, i.e., by splitting the Poisson matrices into several anti-symmetric parts. In section 5, two numerical experiments are done to validate the code, especially energy conservation is demonstrated.  Finally, we conclude this paper.

\section{Laser plasma models with spin effects}
\label{section_reduced_model}
In this section, we introduce the reduced fully relativistic laser-plasma models with spin effects, which are derived based on the conservation of canonical momentum of particles in one and two dimensional case from the three dimensional spin Vlasov--Maxwell model~\cite{LPFEEC, LPS}~(see also in Appendix~\ref{sec:threed}). 
\subsection{One dimensional case}
We assume that an electromagnetic wave is propagating in the longitudinal $x$ direction and that all fields depend spatially on $x$ only.
Choosing the Coulomb gauge $\nabla \cdot \mathbf A =0$, then the vector potential  ${\mathbf A}$ can be denoted as 
${\mathbf A} = { (0, A_y, A_z) = }(0, {\mathbf A}_\perp$). Using ${\mathbf E} = -\nabla \phi - \partial_t {\mathbf A}$, we then obtain with ${\mathbf E} = (E_x, E_y, E_z) =(E_x, {\mathbf E}_\perp)$:
$
{\mathbf E}_\perp = - \partial_t {\mathbf A}_\perp \mbox{ and } E_x=-\partial_x \phi.
$
As for the distribution function, as now the system only depends on $x$ in space, we know that the $y, z$ components of canonical momentum are constants for each particle, i.e., $p_y + A_y, p_z + A_z$ are both constants. When the constants are 0, we get the following one dimensional reduced model, 
The longitudinal variable $p_x$ will be simply denoted by $p$ for convienience.
\begin{equation}
\label{eq:reduced}
\begin{aligned}
&\frac{\partial f}{\partial t} + \frac{p}{\gamma} \frac{\partial f}{\partial x} + \Big[ E_x  - \frac{{\mathbf A}_\perp}{\gamma} \cdot \frac{\partial {\mathbf A}_\perp}{\partial x} + \textgoth{h} \nabla({\mathbf s} \cdot {\mathbf B}) \Big] \frac{\partial f}{\partial p} + (\mathbf{s}\times \mathbf{B} )\cdot \frac{\partial f}{\partial {\mathbf s}} = 0, \\
&\frac{\partial E_x}{\partial t} = -\int_{\mathbb{R}^4} \frac{p}{\gamma} f  \mathrm{d}p\mathrm{d}\mathrm{\mathbf s},\\
&\frac{\partial E_y}{\partial t} = - \frac{\partial^2 A_y}{\partial x^2} + A_y \int_{\mathbb{R}^4}  \frac{f}{\gamma}  \mathrm{d}p\mathrm{d}\mathrm{\mathbf s} + \textgoth{h}\int_{\mathbb{R}^4} s_z \frac{\partial f}{\partial x}\mathrm{d}p\mathrm{d}\mathrm{\mathbf s},\\
&\frac{\partial E_z}{\partial t} = - \frac{\partial^2 A_z}{\partial x^2} + A_z \int_{\mathbb{R}^4}  \frac{f}{\gamma}  \mathrm{d}p\mathrm{d}\mathrm{\mathbf s} - \textgoth{h} \int_{\mathbb{R}^4} s_y \frac{\partial f}{\partial x}\mathrm{d}p\mathrm{d}\mathrm{\mathbf s},\\
& \frac{\partial {\mathbf A}_\perp}{\partial t} = - {\mathbf E}_\perp,\\
&\frac{\partial E_x}{\partial x} = \int_{\mathbb{R}^4} f \mathrm{d}p\mathrm{d}\mathrm{\mathbf s} - 1, \quad (\text{Poisson equation}),
\end{aligned}
\end{equation}
where $\gamma = \sqrt{1 + p^2 + |{\mathbf A}_\perp|^2}$ is the relativistic factor,  $\textgoth{h}$ is the normalized Planck constant, and $\mathbf{B} = \nabla \times \mathbf{ A} = \left(0, -\frac{\partial A_z}{\partial x}, \frac{\partial A_y}{\partial x} \right)^\top$. We can see that $\gamma $ depends on both $x$ and $p$, which brings some difficulties for energy conservation.
This reduced spin Vlasov--Maxwell system possesses a non-canonical Poisson structure~\cite{LPFEEC}.
For any two functionals $\mathcal{F}$ and $\mathcal{G}$ depending on the unknowns
$f, {\mathbf E}$, and ${\mathbf A}_\perp$, the Poisson bracket  is
\begin{equation}\label{eq:poisson}
\begin{aligned}
\{ \mathcal{F}, \mathcal{G}\} &= \int_{\mathbb{R}^5} f\left[\frac{\delta\mathcal{F}}{\delta f},\frac{\delta\mathcal{G}}
{\delta f}\right]_{{xp}}\mathrm{d}{x}\mathrm{d}{p}\mathrm{d}\mathbf{s}
  +\int_{\mathbb{R}^5} \left(\frac{\delta\mathcal{F}}{\delta {E_x}}\frac{\partial f}{\partial {p}}\frac{\delta\mathcal{G}}{\delta f}-\frac{\delta\mathcal{G}}{\delta {E_x}}\frac{\partial f}{\partial {p}}\frac{\delta\mathcal{F}}{\delta f}\right)\mathrm{d}{x}\mathrm{d}{p}\mathrm{d}\mathbf{s}\\
&\hspace{-1.2cm}+ \int_{\mathbb{R}} \left( \frac{\delta \mathcal{G}}{\delta {\mathbf A}_\perp} \cdot  \frac{\delta \mathcal{F}}{\delta {\mathbf E}_\perp} - \frac{\delta \mathcal{F}}{\delta {\mathbf A}_\perp} \cdot  \frac{\delta \mathcal{G}}{\delta {\mathbf E}_\perp}\right) \mathrm{d}x + \frac{1}{\textgoth{h}}\int_{\mathbb{R}^5} f {\mathbf s}\cdot \left( \frac{\partial}{\partial {\mathbf s}}\frac{\delta \mathcal{F}}{\delta {f}} \times \frac{\partial}{\partial {\mathbf s}}\frac{\delta \mathcal{G}}{\delta {f}} \right) \mathrm{d}x\mathrm{d}p\mathrm{d}{\mathbf s},
\end{aligned}
\end{equation}
and the Hamiltonian functional, which is the sum of kinetic, electric, magnetic and Zeeman (spin-dependent) energies, is 
\begin{equation}
\begin{aligned}
\mathcal{H}(f, {\mathbf E}, {\mathbf A}_\perp) &= \int_{\mathbb{R}^5}  (\sqrt{1+p^2+|{\mathbf A}_\perp|^2} - 1)f \mathrm{d}x\mathrm{d}p\mathrm{d}{\mathbf s} + \frac{1}{2}\int_{\mathbb{R}}  |{\mathbf E}|^2 \mathrm{d}x \\
& +   \frac{1}{2}\int_{\mathbb{R}}  \left|\frac{\partial {\mathbf A}_\perp}{\partial x}\right|^2 \mathrm{d}x + \textgoth{h}\int_{\mathbb{R}^5}  \left(  s_y \frac{\partial A_z}{\partial x} {-} s_z \frac{\partial A_y}{\partial x} \right) f \mathrm{d}x\mathrm{d}p\mathrm{d}{\mathbf s}.
\end{aligned}
\label{total_energy}
\end{equation}
Then the reduced spin Vlasov-Maxwell system of equations~(\ref{eq:reduced}) can be reformulated as
$$
\frac{\partial \mathcal{Z}}{\partial t} = \{ \mathcal{Z}, \mathcal{H} \},
$$
where $\mathcal{Z}=(f, E_x,E_y,E_z,A_y,A_z)$. { In this work, periodic boundary condition for $x$ in a finite domain and vanishing boundary conditions for $p \in \mathbb{R}$ and ${\mathbf s} \in \mathbb{R}^3$ are considered. Initial condition is ${\mathcal Z}(t=0)={\mathcal Z}_0=(f_0, {\mathbf E}_{0}, {\mathbf A}_{\perp, 0})$}. 

\subsection{Two dimensional case}
Similar to the one dimensional reduction, we assume an electromagnetic wave propagating in the longitudinal $x_1, x_2$ direction and assuming that the system depend on $x_1, x_2$ only in space. As for the  distribution function, we assume $p_z = - A_z$. Combined with two dimensional reduced Maxwell's equations, we have the following two dimensional reduced model.
\begin{equation}\label{eq:555}
\begin{aligned}
&\frac{\partial f}{\partial t} + \frac{\bf p}{\tilde{\gamma}}\cdot \frac{\partial f}{\partial {\bf x}}+  \cdot \frac{\partial f}{\partial {\mathbf p}} + \left({\mathbf E}_{\mathbf xy} + \tilde{\mathbf F} + \textgoth{h} \nabla{ ({\mathbf s}\cdot{\mathbf B})}  \right) \frac{\partial f}{\partial {\mathbf p}}  + (\mathbf{s}\times \mathbf{B} )\cdot \frac{\partial f}{\partial {\mathbf s}} =0,\\
&{\mathbf B} = \left( \frac{\partial A_z}{\partial x_2}, -\frac{\partial A_z}{\partial x_1}, B_z\right)^\top ,\quad \tilde{\mathbf F} = \left(\frac{p_y B_z+A_zB_y}{\tilde{\gamma}},  - \frac{p_x B_z+A_zB_x}{\tilde{\gamma}}  \right)^\top\\
& \frac{\partial E_x}{\partial t} = \frac{\partial B_z}{\partial y} - \int \frac{p_x}{\tilde{\gamma}}f \mathrm{d}{\bf p} \mathrm{d}{\bf x} + \textgoth{h} \int s_3 \frac{\partial f}{\partial x_2} \mathrm{d}{\mathbf p} \mathrm{d}{\mathbf s},\\
& \frac{\partial E_y}{\partial t} = -\frac{\partial B_z}{\partial x} - \int \frac{p_y}{\tilde{\gamma}}f \mathrm{d}{\bf p} \mathrm{d}{\bf x} - \textgoth{h} \int s_3 \frac{\partial f}{\partial x_1} \mathrm{d}{\mathbf p} \mathrm{d}{\mathbf s},\\
& \frac{\partial B_z}{\partial t} = \frac{\partial E_x}{\partial y} - \frac{\partial E_y}{\partial x},\\
& \frac{\partial A_z}{\partial t} = - E_z,\\
& \frac{\partial E_z}{\partial t} = -\frac{\partial^2 A_z}{\partial x^2} - \frac{\partial^2 A_z}{\partial y^2} + \int \frac{A_z}{\tilde{\gamma}}f \mathrm{d}{\bf p} \mathrm{d}{\bf x} + \textgoth{h} \int \left( s_2 \frac{\partial f}{\partial x_1} - s_1 \frac{\partial f}{\partial x_2} \right)\mathrm{d}{\mathbf p} \mathrm{d}{\mathbf s},\\
& \nabla_{\mathbf x} \cdot {\mathbf E}_{xy} = \int f  \mathrm{d}{\mathbf p} \mathrm{d}{\mathbf s} - 1, \quad (\text{Poisson equation}),\\
& \tilde{\gamma} = \sqrt{1+|{\mathbf p}|^2 + A_z^2}, 
\end{aligned}
\end{equation}
where ${\mathbf x} = (x_1, x_2)^\top, {\mathbf p} = (p_x, p_y)^\top$, ${\mathbf s} \in \mathbb{R}^3$, and  ${\mathbf E}_{xy} = ({\mathbf E}_x^\top, {\mathbf E}_y^\top)^\top$.
For the above model, we for the first time propose its Poisson bracket as
\begin{equation}
\label{eq:2Dbra}
\begin{aligned}
&\{\mathcal{F},\mathcal{G} \}(f({\bf x}, {\bf p}), A_z, B_z, {\bf E})=
 \int f \left[\frac{\delta \mathcal{F}}{\delta f}, \frac{\delta \mathcal{G}}{\delta f}\right]_{\bf{xp}}\mathrm{d}{\bf x} \mathrm{d}{\bf p} +\int \frac{\delta \mathcal{G}}{\delta B_z}  \nabla \times  \frac{\delta \mathcal{F}}{\delta {\bf E}_{xy}} - \frac{\delta \mathcal{F}}{\delta B_z} \nabla \times  \frac{\delta \mathcal{G}}{\delta {\bf E}_{xy}}\mathrm{d}{\bf x}\mathrm{d}{\bf p}  \\
&+ \int \left(\frac{\delta \mathcal{F}}{\delta E_z}\frac{\delta \mathcal{G}}{\delta A_z} - \frac{\delta \mathcal{G}}{\delta E_z}\frac{\delta \mathcal{F}}{\delta A_z}\right) \mathrm{d}{\bf x} + \int f \left(  \frac{\partial }{\partial {\mathbf p}}\frac{\delta \mathcal{F}}{\delta f} \cdot \frac{\delta \mathcal{G}}{\delta {\bf E}_{xy}}  - \frac{\partial }{\partial {\mathbf p}}\frac{\delta \mathcal{G}}{\delta f} \cdot \frac{\delta \mathcal{F}}{\delta {\bf E}_{xy}}\right) \mathrm{d}{\bf x}\mathrm{d}{\bf p}\\
& + \int fB_z\left(\frac{\partial}{\partial p_x}\frac{\delta \mathcal{F}}{\delta f} \frac{\partial}{\partial p_y}\frac{\delta \mathcal{G}}{\delta f} -  \frac{\partial}{\partial p_y}\frac{\delta \mathcal{F}}{\delta f} \frac{\partial}{\partial p_x}\frac{\delta \mathcal{G}}{\delta f} \right)\mathrm{d}{\bf x}\mathrm{d}{\bf p} + \frac{1}{\textgoth{h}}\int_{\mathbb{R}^5} f {\mathbf s}\cdot \left( \frac{\partial}{\partial {\mathbf s}}\frac{\delta \mathcal{F}}{\delta {f}} \times \frac{\partial}{\partial {\mathbf s}}\frac{\delta \mathcal{G}}{\delta {f}} \right)  \mathrm{d}{\mathbf s} \mathrm{d}{\bf p} \mathrm{d}{\bf x}.
\end{aligned}
\end{equation}
With the following Hamiltonian,
\begin{align*}
\mathcal{H} &= \int \left(\sqrt{1+|{\mathbf p}|^2 + A_z^2} -1\right) f  \mathrm{d}{\mathbf s} \mathrm{d}{\bf p} \mathrm{d}{\bf x} + \textgoth{h}  \int {\mathbf s} \cdot {\mathbf B} f \mathrm{d}{\mathbf s} \mathrm{d}{\mathbf p}\mathrm{d}{\mathbf x},\\
& + \frac{1}{2} \int |{\bf E}|^2  \mathrm{d}{\bf x}  + \frac{1}{2}\int |\nabla^\top A_z|^2   \mathrm{d}{\bf x}+ \frac{1}{2} \int {B_z}^2  \mathrm{d}{\bf x},
\end{align*}
the above 2D reduced model could be written as 
$$
\frac{\partial \mathcal{Z}}{\partial t} = \{ \mathcal{Z}, \mathcal{H} \},
$$
where $\mathcal{Z}=(f, {\mathbf E}, A_z, B_z)$. In the above we use the following operators
$$
\nabla f = (\partial_{x_1}f, \partial_{x_2}f)^\top, \nabla^\top f = (\partial_{x_2}f, -\partial_{x_1}f)^\top, \nabla \times {\mathbf f} =  \partial_{x_1}f_2 - \partial_{x_2}f_1.
$$
Similar to one dimensional reduced model, {periodic boundary condition for ${\mathbf x}$ in a finite domain and vanishing boundary conditions for ${\mathbf p} \in \mathbb{R}^2$ and ${\mathbf s} \in \mathbb{R}^3$ are considered. Initial condition is ${\mathcal Z}(t=0)={\mathcal Z}_0=(f_0, {\mathbf E}_{0}, {A}_{z0}, B_{z0})$}.

\section{Semi-discretization}
\label{numeric}
In this section, we introduce the phase-space discretizations for the above two reduced models briefly in the framework of finite element exterior calculus~\cite{FEEC} and particle-in-cell method.
\subsection{One dimensional case}
Following \cite{GEMPIC, LPFEEC, GEMPIC, STRUPHY}, we discretize the components of the electromagnetic fields differently, and consider $E_x, B_y, B_z$ as 1-forms and $E_y, E_z, A_y, A_z$ as 0-forms, which are discretized in finite element spaces $V_0 \subset H^1$ and $V_1 \subset L^2$ respectively. There exists a commuting diagram~\eqref{diagram}
for the involved functional spaces in one spatial dimension, between continuous spaces in the upper line and
discrete subspaces in the lower line. The projectors
$\Pi_0$ and $\Pi_1$ must be constructed carefully in order to assure the diagram to be commuting, such as the quasi-inter/histopolation detailed in~\cite{STRUPHY}.
\begin{equation}
\large
\begin{aligned}
\label{diagram}
 \xymatrix{
    H^{1} \ar[rr]^{\frac{\mathrm{d}}{\mathrm{d}x}} \ar[d]_{\Pi_0} & & L^2  \ar[d]^{\Pi_1} \\
      V_0  \ar[rr]^{\frac{\mathrm{d}}{\mathrm{d}x}} && V_1
    }
\end{aligned}
\end{equation}
The spatial domain $[0, L]$ is discretized by a uniform grid
$$
x_j = j\Delta x, \;\; \Delta x = L/M, \;\; j = 0, \cdots, M-1.
$$
In this paper, we choose the B-splines~\cite{buffa} of order $k, k-1, k \ge 1$ on the above uniform meshes as the basis functions for $V_0$, and $V_1$ with periodic boundary condition, which are denoted as 
$\{ \Lambda^0_j(x)\}_{j=0,\cdots,N_0-1} $, and $\{ \Lambda^1_j(x)\}_{j=0,\cdots,N_1-1} $, i.e.,  
$$\Lambda^0_j(x) = N^k_j(x), \quad \Lambda^1_j(x) = N^{k-1}_j(x), 0 \le j < M(=N_0=N_1),$$
where $N_i^k$ is the B-splines of degree $k$ with the support of $[x_i, \cdots, x_{i+k+1})$.
The important relation between  $\Lambda^1$ and $\Lambda^0$:
$
\frac{\mathrm{d}}{\mathrm{d}x} \Lambda^0_j(x) = \frac{1}{\Delta x}\left(\Lambda^1_j(x) - \Lambda^1_{j+1}(x)  \right)
$
 can be reformulated as
$$
{ \frac{\mathrm{d}}{\mathrm{d}x}(\Lambda^0_0, \cdots, \Lambda^0_{N_0-1})(x)  = (\Lambda^1_0, \cdots, \Lambda^1_{N_1-1})(x) \mathbb{G}},
$$
where the size of matrix $\mathbb{G}$ is $ N_1 \times N_0$.
The approximations of electric field and magnetic potential components can be written as
\begin{equation}
\label{fem_e}
E_{x,h}(t, x) = \sum_{j=0}^{N_1-1} e_{x, j}(t) \Lambda^1_j(x), \; E_{y,h}(t, x) = \sum_{j=0}^{N_0-1} e_{y, j}(t) \Lambda^0_j(x),
\; E_{z,h}(t, x) = \sum_{j=0}^{N_0-1} e_{z, j}(t) \Lambda^0_j(x),
\end{equation}
\begin{equation}
\label{fem_a}
A_{y,h}(t, x) = \sum_{j=0}^{N_0-1} { a_{y, j}(t)} \Lambda^0_j(x),  \;
A_{z,h}(t, x) = \sum_{j=0}^{N_0-1} { a_{z, j}(t)} \Lambda^0_j(x).
\end{equation}
The distribution function $f(t, x, p, \mathbf{s})$ is discretized as the sum of finite number of particles with constant weights, i.e., 
\begin{equation}
\label{f_dirac}
f(t, x, p, \mathbf{s}) \approx f_h(t, x, p, \mathbf{s}) = \sum_{a=1}^{N_p} \omega_a \delta(x-x_a(t))\delta(p-p_a(t))\delta(\mathbf{s}-\mathbf{s}_a(t)),
\end{equation}
where $N_p$ is the total particle number, $\omega_a$, $x_a$, $p_a$, and $\mathbf{s}_a$ denote the weight,  the position, the momentum (velocity), and the
spin co-ordinates of $a$-th particle, respectively, $1 \le a \le N_p$.

By discretizing the Poisson bracket using discrete functional derivatives as in~\cite{LPFEEC}, we have the following discrete Poisson bracket.
\begin{equation}\label{eq:dis}
\{ F, G \} = \left( \nabla_{\mathbf u}F\right)^{\top} \mathbb{J}({\mathbf u})\nabla_{\mathbf u}G,
\end{equation}
where ${\mathbf u} = ({\mathbf X}, {\mathbf P}, {\mathbf S}, {\mathbf e}_x, {\mathbf e}_y, {\mathbf e}_z, {\mathbf a}_y, {\mathbf a}_z)^{\top}$
and  the matrix $\mathbb{J}({\mathbf u})$ is defined by
\begin{align}\label{eq:dispoisson}
\begin{aligned}
& \mathbb{J}({\mathbf{u}}) \!=\!\!
\left(
\begin{matrix}
    {\mathbf{0}} \!\! & {\mathbb{W}}^{-1}  &{\mathbf{0}} & {\mathbf{0}} \!\! & {\mathbf{0}} & {\mathbf{0}}& {\mathbf{0}}  & {\mathbf{0}}\\
     -{\mathbb{W}}^{-1}  \!\! &{\mathbf{0}}   & {\mathbf{0}} & \mathbb{\Lambda}_1({\mathbf {X}})\mathbb{M}_1^{-1} \!\! & {\mathbf{0}} & {\mathbf{0}} & {\mathbf{0}} & {\mathbf{0}}\\
      {\mathbf{0}}  \!\! & {\mathbf{0}}  & \frac{1}{\textgoth{h}}\mathbb{S} & {\mathbf{0}} \!\!& {\mathbf{0}}  & {\mathbf{0}} & {\mathbf{0}} & {\mathbf{0}} \\
       {\mathbf{0}}  \!\! & -{\mathbb{M}}_1^{-1}\mathbb{\Lambda}_1({\mathbf{X}})^{\top}  & {\mathbf{0}}& {\mathbf{0}}  \!\!& {\mathbf{0}}& {\mathbf{0}}& {\mathbf{0}}& {\mathbf{0}}\\
       {\mathbf{0}}\!\! & {\mathbf{0}}& {\mathbf{0}}& {\mathbf{0}}\!\!& {\mathbf{0}}& {\mathbf{0}} &{\mathbb{M}}_0^{-1}&{\mathbf{0}}\\
       {\mathbf{0}}\!\! &{\mathbf{0}}&{\mathbf{0}}&{\mathbf{0}}\!\!&{\mathbf{0}}&{\mathbf{0}} &{\mathbf{0}}& {\mathbb{M}}_0^{-1}\\
       {\mathbf{0}}\!\! &{\mathbf{0}}&{\mathbf{0}}&{\mathbf{0}}\!\!&-{\mathbb{M}}_0^{-1} &{\mathbf{0}}&{\mathbf{0}}&{\mathbf{0}}\\
       {\mathbf{0}}\!\! &{\mathbf{0}}&{\mathbf{0}}&{\mathbf{0}}\!\!&{\mathbf{0}} &-{\mathbb{M}}_0^{-1}&{\mathbf{0}}&{\mathbf{0}}
      \end{matrix}
\right).
\end{aligned}
\end{align}
In the above, ${\mathbf{X}}=(x_a), {\mathbf{P}}=(p_a), {\mathbf{S}}=({\mathbf s}_a)$ denote three vectors of sizes $N_p, N_p, 3N_p$ storing the positions, velocities, and spin values of all particles. ${\mathbf S}_x$, ${\mathbf S}_y$, and ${\mathbf S}_z$ denote three long vectors storing the $x, y, z$ component of spin variable of all particles. ${\mathbf{e}}_x = (e_{x,i})$, ${\mathbf{e}}_y = (e_{y,i})$, ${\mathbf{e}}_z = (e_{z,i})$, ${\mathbf{a}}_y = (a_{y,i})$, and ${\mathbf{a}}_x = (a_{z,i})$ denote the degrees of freedoms of fields. $\mathbb{\Lambda}_i({\mathbf X})$ is a matrix of size $N_p \times N_i$ storing the values of basis functions of $V_i, i = 0, 1$ evaluated at all the particle positions. ${\Lambda}_i(x_a)$ means a vector storing the values of basis functions of $V_i, i = 0, 1$ at $a$-th particle position. $\mathbb{M}_i$ is the mass matrix of finite element space $V_i, i = 0, 1$.
Finally, we introduce the weight matrix $\mathbb{W} =
\text{diag}(\omega_1, \cdots, \omega_{N_p})\in {\cal M}_{N_p, N_p}(\mathbb{R})$, and $\mathbb{S} = \text{diag}(\mathbb{S}_1, \cdots, \mathbb{S}_{N_p})\in {\cal M}_{3N_p, 3N_p}(\mathbb{R})$, where
$$
\mathbb{S}_a = \frac{1}{\omega_a}
\left(
\begin{matrix}
     0 & s_{a,z}  & -s_{a,y} \\
      -s_{a,z}  & 0  & s_{a,x} \\
      s_{a,y}  & -s_{a,x}  & 0
      \end{matrix}
\right)\in {\cal M}_{3, 3}(\mathbb{R}), \quad 1 \le a \le N_p.
$$
%
Using the notations, discrete Hamiltonian can be written more compactly as
\begin{equation}
\label{discrete_H}
\begin{aligned}
H({\mathbf u})&\displaystyle= \sum_{a=1}^{N_p}\omega_a ( \sqrt{1 + p_a^2 + |{\mathbf A}_{\perp, h}(x_a)|^2} -1 )  + {\textgoth{h}} \, {\mathbf a}_z^{\top}\mathbb{G}^{\top}\mathbb{\Lambda}_1({\mathbf X})^{\top}\mathbb{W}{\mathbf S}_y
- {\textgoth{h}}\, {\mathbf a}_y^{\top}\mathbb{G}^{\top}\mathbb{\Lambda}_1({\mathbf X})^{\top}\mathbb{W}{\mathbf S}_z,\\
&\displaystyle + \frac{1}{2}{\mathbf e}_x^{\top}\mathbb{M}_1{\mathbf e}_x +\frac{1}{2}{\mathbf e}_y^{\top}\mathbb{M}_0{\mathbf e}_y +  \frac{1}{2}{\mathbf e}_z^{\top}\mathbb{M}_0{\mathbf e}_z +  \frac{1}{2}{\mathbf a}_y^{\top} \mathbb{G}^{\top}\mathbb{M}_1\mathbb{G}{\mathbf a}_y+  \frac{1}{2}{\mathbf a}_z^{\top} \mathbb{G}^{\top}\mathbb{M}_1\mathbb{G}{\mathbf a}_z.
\end{aligned}
\end{equation}
From the discrete Poisson bracket \eqref{eq:dis}-\eqref{eq:dispoisson} and the discrete Hamiltonian \eqref{discrete_H},
the equations of motion then read as
\begin{equation}
\label{ode_1}
\dot{\mathbf u} = \left\{ {\mathbf u}, H\right\} = \mathbb{J}({\mathbf u})\nabla_{\mathbf u}H, \;\;\; {\mathbf u}(t=0) = {\mathbf u}_0.
\end{equation}

\subsection{Two dimensional case}
In the two dimensional case, we regards $A_z, E_z$ are 0-forms, ${\mathbf E}_{xy}$ as a 1-form, $B_z$ as a 2-form, and corresponding finite element spaces make the following diagram commute with suitable projectors,
\begin{equation}
\large
\begin{aligned}
\label{diagram2d}
 \xymatrix{
    H^{1} \ar[rr]^{\nabla} \ar[d]_{\Pi_0} & & H(\text{curl})  \ar[d]^{\Pi_1} \ar[rr]^{\nabla \times}  && L^2  \ar[d]^{\Pi_2} \\
      V_0  \ar[rr]^{\nabla} && V_1   \ar[rr]^{\nabla \times}   && V_2
    }
\end{aligned}
\end{equation}
In the following, we describe the discretization with a slight abuse of notation with one dimensional notations. 
We assume a uniform grid on spatial domain $[0, L_1] \times [0, L_2]$ with 
$$
x_{i, j} = j\Delta x_i, \;\; \Delta x_i = L_i/M_i, \;\; j = 0, \cdots, M_i-1, \quad i = 1, 2.
$$
The basis functions of $V_i, i = 0, 1, 2$ are the tensor products of B-splines, i.e., 
\begin{equation}
\begin{aligned}
&V_0:= \text{span} \{ \Lambda^0_i \ |\  0 \le i < N_{V_0}=M_1M_2 \},  \Lambda^0_i({\mathbf x}) := N^{k_1}_{i_1}(x_1)N^{k_2}_{i_2}(x_2), \ i = i_1M_2 + i_2,\\
 &V_1 := \text{span}  \left\{
\left(\begin{matrix}
  \Lambda^1_{1,i} \\
    0
\end{matrix} \right),
\left(\begin{matrix}
  0\\
     \Lambda^1_{2,i}
\end{matrix} \right)
\Bigg| 
\begin{matrix}
  0 \le i < N_{V_1, 1} =M_1M_2 \\
    0 \le i < N_{V_1, 2} = M_1M_2
\end{matrix} 
\right\},\\
& \Lambda^1_{1,i} := N^{k_1-1}_{i_1}(x_1)N^{k_2}_{i_2}(x_2), \quad \Lambda^1_{2,i} := N^{k_1}_{i_1}(x_1)N^{k_2-1}_{i_2}(x_2),\ i = i_1M_2 + i_2,\ N_{V_1} = N_{V_1,1} + N_{V_1,2},\\
&V_2:=  \text{span} \{ \Lambda^2_i \ |\  0 \le i < N_{V_2} = M_1M_2\}, \Lambda^2_i({\mathbf x}) := N^{k_1-1}_{i_1}(x_1)N^{k_2-1}_{i_2}(x_2), \ i = i_1M_2 + i_2.
\end{aligned}
\end{equation}
We also introduce another finite element space denoted by $V_1^*$ (where $(B_x, B_y)^\top$ is discretized) with following basis functions
\begin{equation}
\begin{aligned}
 &V_1^* := \text{span}  \left\{
\left(\begin{matrix}
  \Lambda^*_{2,i} \\
    0
\end{matrix} \right),
\left(\begin{matrix}
  0\\
    - \Lambda^*_{1,i}
\end{matrix} \right)
\Bigg| 
\begin{matrix}
  0 \le i < N^{*}_1 = M_1M_2\\
    0 \le i < N^{*}_2 = M_1M_2
\end{matrix} 
\right\},\\
& \Lambda^*_{1,i} := N^{k_1-1}_{i_1}(x_1)N^{k_2}_{i_2}(x_2), \quad \Lambda^*_{2,i} := N^{k_1}_{i_1}(x_1)N^{k_2-1}_{i_2}(x_2),\ i = i_1M_2 + i_2, \ N^* = N^*_1 + N^*_2.
\end{aligned}
\end{equation}
The matrices of linear operator $\nabla$, $\nabla \times $, and $\nabla^\top$ are denoted as $\mathbb{G}$, ${\mathbb{C}}$, and $\mathbb{G}_{*}$ with sizes $N_{V_1} \times N_{V_0}, N_{V_2} \times N_{V_1}, N^* \times N_{V_0}$, respectively.
Mass matrices of $V_0, V_1, V_2, V^*_1$ are denoted as ${\mathbb{M}}_0, {\mathbb{M}}_1, {\mathbb{M}}_2, {\mathbb{M}}_{1,*}$, respectively. 
$\mathbb{\Lambda}_i({\mathbf X})$ is a matrix of size $N_p \times N_i$ storing the values of basis functions of $V_i, i = 0, 1$ evaluated at all the particle positions. $\mathbb{\Lambda}_1({\mathbf X}) \ (\mathbb{\Lambda}_*({\mathbf X}))$ is a matrix of size $2N_p \times N_1 \ (2N_p \times N^*)$ storing the values of basis functions of $V_1 \ V^*_1$ evaluated at all the particle positions. 
$\Lambda_i({\mathbf x}_a)$ denotes a vector of length of $N_{V_i}$ storing the values of all the basis functions of $V_{i},  i = 0, 2$ at $a$-th particle positions.   $\Lambda_1({\mathbf x}_a)$ ($\Lambda_*({\mathbf x}_a)$) denotes a matrix of size of $N_{V_1} \times 2 \ (N^* \times 2)$ storing the values of all the basis functions of $V_{1} \ (V_1^*)$ at $a$-th particle position. Distribution function is discretized as the sum of $N_p$ particles with constant weights as~\eqref{f_dirac}.
By discretizing functional derivatives~(see in appendix~\ref{sec:dis2d}) in~\eqref{eq:2Dbra}, we get the following discrete Poisson bracket
\begin{equation}\label{eq:dis2D}
\{ F, G \} = \left( \nabla_{\mathbf u}F\right)^{\top} \mathbb{J}({\mathbf u})\nabla_{\mathbf u}G,
\end{equation}
where ${\mathbf u} = ({\mathbf X}, {\mathbf P}, {\mathbf S}, {\mathbf e}_{xy}, {\mathbf b}_{z}, {\mathbf e}_z, {\mathbf a}_z)^{\top}$
and  the matrix $\mathbb{J}({\mathbf u})$ is defined by
\begin{align}\label{eq:dispoisson2d}
\begin{aligned}
& \mathbb{J}({\mathbf{u}}) \!=\!\!
\left(
\begin{matrix}
    {\mathbf{0}} \!\! & {{\mathbb{W}}}^{-1}  &{\mathbf{0}} & {\mathbf{0}} \!\! & {\mathbf{0}} & {\mathbf{0}}& {\mathbf{0}}  \\
     -{{\mathbb{W}}}^{-1}  \!\! & \mathbb{S}^p  & {\mathbf{0}} & \mathbb{\Lambda}_1({\mathbf {X}}){\mathbb{M}}_1^{-1} \!\! & {\mathbf{0}} & {\mathbf{0}} & {\mathbf{0}} \\
      {\mathbf{0}}  \!\! & {\mathbf{0}}  & \frac{1}{\textgoth{h}} {\mathbb{S}} & {\mathbf{0}} \!\!& {\mathbf{0}}  & {\mathbf{0}} & {\mathbf{0}} \\
       {\mathbf{0}}  \!\! & -{{\mathbb{M}}}_1^{-1}\mathbb{\Lambda}_1({\mathbf{X}})^{\top}  & {\mathbf{0}}& {\mathbf{0}}  \!\!& {\mathbb{M}}_1^{-1}{\mathbb{C}}^\top& {\mathbf{0}}& {\mathbf{0}}\\
       {\mathbf{0}}\!\! & {\mathbf{0}}& {\mathbf{0}}& -{\mathbb{C}} {{\mathbb{M}}}_1^{-1} \!\!& {\mathbf{0}}& {\mathbf{0}} &  {\mathbf{0}}\\
       {\mathbf{0}}\!\! &{\mathbf{0}}&{\mathbf{0}}&{\mathbf{0}}\!\!&{\mathbf{0}}&{\mathbf{0}} &{{\mathbb{M}}}_0^{-1}\\
       {\mathbf{0}}\!\! &{\mathbf{0}}&{\mathbf{0}}&{\mathbf{0}}\!\!&  {\mathbf{0}} &-{{\mathbb{M}}}_0^{-1}&{\mathbf{0}}
               \end{matrix}
\right),
\end{aligned}
\end{align}
where $ {\mathbf e}_{xy}$ is a vector storing the finite element coefficients of $(E_{x,h}, E_{y,h})^\top$ (a one form), ${\mathbf{X}}=({\mathbf x}_a), {\mathbf{P}}=({\mathbf p}_a), {\mathbf{S}}=({\mathbf s}_a)$ denote three vectors of sizes $2N_p, 2N_p, 3N_p$ storing the positions, velocities, and spin values of all particles.
and
\begin{align}
\begin{aligned}
&\mathbb{S}^p_a = \frac{1}{\omega_a}
\left(
\begin{matrix}
     0 & B_{z,h}({\mathbf x}_a)   \\
  -B_{z,h}({\mathbf x}_a)   & 0 
      \end{matrix}
\right)\in {\cal M}_{2, 2}(\mathbb{R}), \quad \mathbb{S}^p = \text{diag}(\mathbb{S}^p_1, \cdots, \mathbb{S}^p_{N_p})\in {\cal M}_{2N_p, 2N_p}(\mathbb{R}).
\end{aligned}
\end{align}
Discrete Hamiltonian is 
\begin{equation}
\label{discrete_H2d}
\begin{aligned}
H({\mathbf u})&\displaystyle= \sum_{a=1}^{N_p}\omega_a ( \sqrt{1 + |{\mathbf p}_a|^2 + |{A}_{z, h}({\mathbf x}_a)|^2} -1 )  + \frac{1}{2} {\mathbf b}_z^\top {\mathbb{M}}_2 {\mathbf b}_z +   \frac{1}{2} {\mathbf e}_{xy}^\top {\mathbb{M}}_1 {\mathbf e}_{xy}  +   \frac{1}{2} {\mathbf e}_{z}^\top {\mathbb{M}}_0 {\mathbf e}_{z}   \\
& + \frac{1}{2} {\mathbf a}_z^\top   \mathbb{G}_{*}^\top {\mathbb{M}}_{1, *} \mathbb{G}_{*} {\mathbf a}_z + {\textgoth{h}}\, {\mathbf a}_z^\top \mathbb{G}^\top_{*}  \mathbb{\Lambda}_{1, *}({\mathbf X})^{\top}{\mathbb{W}}{\mathbf S}_{xy} +  {\textgoth{h}}\, {\mathbf b}_z^{\top} \mathbb{\Lambda}_2({\mathbf X})^{\top}{\mathbb{W}}{\mathbf S}_z,
\end{aligned}
\end{equation}
where ${\mathbf S}_{xy} = ({\mathbf S}_x^\top, {\mathbf S}_y^\top)^\top$, 
with which we obtain a time-continuous Poisson system 
\begin{equation}
\label{ode_2}
\dot{\mathbf u} = \left\{ {\mathbf u}, H\right\} = \mathbb{J}({\mathbf u})\nabla_{\mathbf u}H, \;\;\; {\mathbf u}(t=0) = {\mathbf u}_0.
\end{equation}

\section{Time discretization}
In this section, as the Hamiltonian splitting method does not give explicitly solvable subsystems, we use Poisson splitting (to split the Poisson matrix) and obtain several subsystems as~\cite{2020en}. By using discrete gradient method proposed in~\cite{DIS, Gonzalez}, energy is conserved by the fully discrete scheme constructed, also discrete Poisson equation is satisfied by the numerical solution. 

\subsection{Discrete gradient method}
For the following form conservative ordinary equations, 
$$
\dot{y} = J({\mathbf y}) \nabla H, \quad J({\mathbf y}) ^\top = - J({\mathbf y}),
$$
$\bar{\nabla}H({\mathbf y}^n, {\mathbf y}^{n+1})$ is called a discrete gradient for the above equations for time step $[t_n, t_{n+1}]$, if 
$$
({\mathbf y}^{n+1} - {\mathbf y}^n)^\top \bar{\nabla}H({\mathbf y}^n, {\mathbf y}^{n+1}) = H({\mathbf y}^{n+1}) -  H({\mathbf y}^{n}). 
$$
Then we obtain the following energy conserving schemes with the help of the discrete gradient,
$$
\frac{{\mathbf y}^{n+1} - {\mathbf y}^{n}   }{\Delta t} = \bar{J}({\mathbf y}^n, {\mathbf y}^{n+1}) \bar{\nabla}H({\mathbf y}^n, {\mathbf y}^{n+1}),
$$
where $\bar{J}({\mathbf y}^n, {\mathbf y}^{n+1})$ is any anti-symmetric approximation of $ J({\mathbf y})$.

\subsection{One dimensional case}
The Poisson matrix~\eqref{eq:dispoisson} is split into following three parts,
\begin{align*}
\begin{aligned}
& \mathbb{J}_1({\mathbf{u}}) \!=\!\!
\left(
\begin{matrix}
    {\mathbf{0}} \!\! & {\mathbb{W}}^{-1}  &{\mathbf{0}} & {\mathbf{0}} \!\! & {\mathbf{0}} & {\mathbf{0}}& {\mathbf{0}}  & {\mathbf{0}}\\
     -{\mathbb{W}}^{-1}  \!\! &{\mathbf{0}}   & {\mathbf{0}} & \mathbb{\Lambda}_1({\mathbf {X}})\mathbb{M}_1^{-1} \!\! & {\mathbf{0}} & {\mathbf{0}} & {\mathbf{0}} & {\mathbf{0}}\\
      {\mathbf{0}}  \!\! & {\mathbf{0}}  &{\mathbf 0} & {\mathbf{0}} \!\!& {\mathbf{0}}  & {\mathbf{0}} & {\mathbf{0}} & {\mathbf{0}} \\
       {\mathbf{0}}  \!\! & -{\mathbb{M}}_1^{-1}\mathbb{\Lambda}_1({\mathbf{X}})^{\top}  & {\mathbf{0}}& {\mathbf{0}}  \!\!& {\mathbf{0}}& {\mathbf{0}}& {\mathbf{0}}& {\mathbf{0}}\\
       {\mathbf{0}}\!\! & {\mathbf{0}}& {\mathbf{0}}& {\mathbf{0}}\!\!& {\mathbf{0}}& {\mathbf{0}} &{\mathbf{0}}&{\mathbf{0}}\\
       {\mathbf{0}}\!\! &{\mathbf{0}}&{\mathbf{0}}&{\mathbf{0}}\!\!&{\mathbf{0}}&{\mathbf{0}} &{\mathbf{0}}&{\mathbf{0}}\\
       {\mathbf{0}}\!\! &{\mathbf{0}}&{\mathbf{0}}&{\mathbf{0}}\!\!&{\mathbf{0}}&{\mathbf{0}}&{\mathbf{0}}&{\mathbf{0}}\\
       {\mathbf{0}}\!\! &{\mathbf{0}}&{\mathbf{0}}&{\mathbf{0}}\!\!&{\mathbf{0}} &{\mathbf{0}}&{\mathbf{0}}&{\mathbf{0}}
      \end{matrix}
\right),
\end{aligned}\\
\begin{aligned}
& \mathbb{J}_2({\mathbf{u}}) \!=\!\!
\left(
\begin{matrix}
    {\mathbf{0}} \!\! & {\mathbf{0}} &{\mathbf{0}} & {\mathbf{0}} \!\! & {\mathbf{0}} & {\mathbf{0}}& {\mathbf{0}}  & {\mathbf{0}}\\
    {\mathbf{0}} \!\! &{\mathbf{0}}   & {\mathbf{0}} &{\mathbf{0}} \!\! & {\mathbf{0}} & {\mathbf{0}} & {\mathbf{0}} & {\mathbf{0}}\\
      {\mathbf{0}}  \!\! & {\mathbf{0}}  & {\mathbf{0}}& {\mathbf{0}} \!\!& {\mathbf{0}}  & {\mathbf{0}} & {\mathbf{0}} & {\mathbf{0}} \\
       {\mathbf{0}}  \!\! &{\mathbf{0}} & {\mathbf{0}}& {\mathbf{0}}  \!\!& {\mathbf{0}}& {\mathbf{0}}& {\mathbf{0}}& {\mathbf{0}}\\
       {\mathbf{0}}\!\! & {\mathbf{0}}& {\mathbf{0}}& {\mathbf{0}}\!\!& {\mathbf{0}}& {\mathbf{0}} &{\mathbb{M}}_0^{-1}&{\mathbf{0}}\\
       {\mathbf{0}}\!\! &{\mathbf{0}}&{\mathbf{0}}&{\mathbf{0}}\!\!&{\mathbf{0}}&{\mathbf{0}} &{\mathbf{0}}& {\mathbb{M}}_0^{-1}\\
       {\mathbf{0}}\!\! &{\mathbf{0}}&{\mathbf{0}}&{\mathbf{0}}\!\!&-{\mathbb{M}}_0^{-1} &{\mathbf{0}}&{\mathbf{0}}&{\mathbf{0}}\\
       {\mathbf{0}}\!\! &{\mathbf{0}}&{\mathbf{0}}&{\mathbf{0}}\!\!&{\mathbf{0}} &-{\mathbb{M}}_0^{-1}&{\mathbf{0}}&{\mathbf{0}}
      \end{matrix}
\right),
\end{aligned}
\begin{aligned}
& \mathbb{J}_3({\mathbf{u}}) \!=\!\!
\left(
\begin{matrix}
    {\mathbf{0}} \!\! &{\mathbf{0}} &{\mathbf{0}} & {\mathbf{0}} \!\! & {\mathbf{0}} & {\mathbf{0}}& {\mathbf{0}}  & {\mathbf{0}}\\
     {\mathbf{0}}\!\! &{\mathbf{0}}   & {\mathbf{0}} &{\mathbf{0}}\!\! & {\mathbf{0}} & {\mathbf{0}} & {\mathbf{0}} & {\mathbf{0}}\\
      {\mathbf{0}}  \!\! & {\mathbf{0}}  & \frac{1}{\textgoth{h}}\mathbb{S} & {\mathbf{0}} \!\!& {\mathbf{0}}  & {\mathbf{0}} & {\mathbf{0}} & {\mathbf{0}} \\
       {\mathbf{0}}  \!\! & {\mathbf{0}} & {\mathbf{0}}& {\mathbf{0}}  \!\!& {\mathbf{0}}& {\mathbf{0}}& {\mathbf{0}}& {\mathbf{0}}\\
       {\mathbf{0}}\!\! & {\mathbf{0}}& {\mathbf{0}}& {\mathbf{0}}\!\!& {\mathbf{0}}& {\mathbf{0}} &{\mathbf{0}}&{\mathbf{0}}\\
       {\mathbf{0}}\!\! &{\mathbf{0}}&{\mathbf{0}}&{\mathbf{0}}\!\!&{\mathbf{0}}&{\mathbf{0}} &{\mathbf{0}}&{\mathbf{0}}\\
       {\mathbf{0}}\!\! &{\mathbf{0}}&{\mathbf{0}}&{\mathbf{0}}\!\!&{\mathbf{0}}&{\mathbf{0}}&{\mathbf{0}}&{\mathbf{0}}\\
       {\mathbf{0}}\!\! &{\mathbf{0}}&{\mathbf{0}}&{\mathbf{0}}\!\!&{\mathbf{0}} &{\mathbf{0}}&{\mathbf{0}}&{\mathbf{0}}
      \end{matrix}
\right),
\end{aligned}
\end{align*}
which correspond to the following three subsystems.

\noindent{\bf Subsystem I}
The first subsystem about variables $x_a, p_a, {\mathbf e}_x, 1 \le a \le N_p$ is 
\begin{equation}
\begin{aligned}
&\dot{x}_a = \frac{p_a}{\sqrt{1+p_a^2 + |{\mathbf A}_\perp(x_a)|^2}},\\
& \dot{p}_a = E_{x, h}(x_a) - \frac{{{\mathbf A}_{\perp, h}}(x_a) \cdot \partial_x {\mathbf A}_{\perp, h}(x_a)} {\sqrt{1+p_a^2 + |{\mathbf A}_\perp(x_a)|^2}} -  {\textgoth{h}} s^a_y  \partial_x^2 A_z(x_a)   +   {\textgoth{h}} s^a_z \partial_x^2 A_y(x_a),\ 1 \le a \le N_p,\\
&\dot{\mathbf e}_x = -\mathbb{M}_1^{-1} \sum_{a=1}^{N_p} {\Lambda}_1(x_a) w_a \frac{p_a}{\sqrt{1+p_a^2 + |{\mathbf A}_\perp(x_a)|^2}},\\
&\dot{\mathbf S} =  {\mathbf 0},\ \dot{\mathbf e}_y = {\mathbf 0},\  \dot{\mathbf e}_z = {\mathbf 0},\ \dot{\mathbf a}_y = {\mathbf 0},\ \dot{\mathbf a}_z = {\mathbf 0}.
\end{aligned}
\end{equation}
For variables $x_a, p_a, {\mathbf e}_x, 1 \le a \le N_p$, we have the following discrete gradients using method in~\cite{Gonzalez}, 
\begin{equation}
\begin{aligned}
&\bar{\nabla}_{x_a} H =  w_a\frac{ \left( {{\mathbf A}^n_{\perp,h}(x_a^{n+1})  + {\mathbf A}^n_{\perp,h}(x_a^n) } \right) \cdot \frac{\left( {{\mathbf A}^n_{\perp,h}(x_a^{n+1})  - {\mathbf A}^n_{\perp,h}(x_a^n) } \right)}{(x_a^{n+1} - x_a^n)}    }{\sqrt{ 1 + (p_a^n)^2 + |{\mathbf A}^n_{\perp,h}(x_a^{n+1})|^2 } + \sqrt{ 1 + (p_a^{n})^2 + |{\mathbf A}^n_{\perp,h}(x_a^{n})|^2 } },\\
& + w_a {\textgoth{h}} s^n_{a,y} \left( \frac{\partial_x A^n_{z,h}(x_a^{n+1}) -  \partial_x A^n_{z,h}(x_a^{n}) }{x_a^{n+1} - x_a^n}  \right) +   {\textgoth{h}} s^n_{a,z} \left( \frac{\partial_x A^n_{y,h}(x_a^{n+1}) -  \partial_x A^n_{y,h}(x_a^{n}) }{x_a^{n+1} - x_a^n}  \right),\\
&\bar{\nabla}_{p_a} H = w_a \frac{ p_a^{n} + p_a^{n+1}}{\sqrt{ 1 + (p_a^n)^2 + |{\mathbf A}^n_{\perp,h}(x_a^{n+1})|^2 } + \sqrt{ 1 + (p_a^{n+1})^2 + |{\mathbf A}^n_{\perp,h}(x_a^{n+1})|^2 }}\\
&\bar{\nabla}_{{\mathbf e}_x} H =\mathbb{M}_1  \frac{{\mathbf e}_x^n + {\mathbf e}_x^{n+1} }{2}.
\end{aligned}
\end{equation}
With the above discrete gradient, we have the following scheme
\begin{equation}
\begin{aligned}
&\dot{ \frac{x_a^{n+1} - x_a^n}{\Delta t} } = \frac{ p_a^{n} + p_a^{n+1}}{\sqrt{ 1 + (p_a^n)^2 + |{\mathbf A}^n_{\perp,h}(x_a^{n+1})|^2 } + \sqrt{ 1 + (p_a^{n+1})^2 + |{\mathbf A}^n_{\perp,h}(x_a^{n+1})|^2 }},\\
& \frac{{ p}_a^{n+1} - p_a^n}{\Delta t} =\frac{1}{\Delta t}\int_{t^n}^{t^{n+1}} { \mathbb{\Lambda}_1({\mathbf X}(\tau))} \mathrm{d}\tau\frac{{\mathbf e}_x^{n+1} + {\mathbf e}_x^n}{2} \\
&- \frac{ \left( {{\mathbf A}^n_{\perp,h}(x_a^{n+1})  + {\mathbf A}^n_{\perp,h}(x_a^n) } \right) \cdot \frac{ {{\mathbf A}^n_{\perp,h}(x_a^{n+1})  - {\mathbf A}^n_{\perp,h}(x_a^n) } }{x_a^{n+1} - x_a^n}    }{\sqrt{ 1 + (p_a^n)^2 + |{\mathbf A}^n_{\perp,h}(x_a^{n+1})|^2 } + \sqrt{ 1 + (p_a^{n})^2 + |{\mathbf A}^n_{\perp,h}(x_a^{n})|^2 } },\\
& -  {\textgoth{h}} s_{a,y}^n \left( \frac{\partial_x A^n_{z,h}(x_a^{n+1}) -  \partial_x A^n_{z,h}(x_a^{n}) }{x_a^{n+1} - x_a^n}  \right) +   {\textgoth{h}} s_{a,z}^n \left( \frac{\partial_x A^n_{y,h}(x_a^{n+1}) -  \partial_x A^n_{y,h}(x_a^{n}) }{x_a^{n+1} - x_a^n}  \right),\\
&\frac{{\mathbf e}_x^{n+1} - {\mathbf e}_x^{n} }{\Delta t} = -\mathbb{M}_1^{-1}\sum_a \frac{1}{\Delta t} \int_{t^n}^{t^{n+1}} {\Lambda}_1(x_a(\tau))^{\top} \mathrm{d}\tau w_a  \\&\frac{ p_a^{n} + p_a^{n+1}}{\sqrt{ 1 + (p_a^n)^2 + |{\mathbf A}^n_{\perp,h}(x_a^{n+1})|^2 } + \sqrt{ 1 + (p_a^{n+1})^2 + |{\mathbf A}^n_{\perp,h}(x_a^{n+1})|^2 }},\\
& {\mathbf S}^{n+1} = {\mathbf S}^{n}  ,\  {\mathbf e}_y^{n+1} = {\mathbf e}_y^{n} , \ {\mathbf e}_z^{n+1}={\mathbf e}_z^{n}  ,\ {\mathbf a}_y^{n+1}  = {\mathbf a}_y^{n} ,\ {\mathbf a}_z^{n+1}  ={\mathbf a}_z^{n},
\end{aligned}
\end{equation}
where the time-continuous trajectory is defined as 
$$
x_a(\tau) = x_a^n + (\tau - t^n) \frac{x_a^{n+1} - x_a^n}{\Delta t}, \quad \tau \in [t^n, t^{n+1}], \quad 1 \le a \le N_p.
$$
\begin{remark}\label{re:dispoisson}
When $x_a^n$ is very close to $x_a^{n+1}$,  $\frac{ {{\mathbf A}^n_{\perp,h}(x_a^{n+1})  - {\mathbf A}^n_{\perp,h}(x_a^n) } }{x_a^{n+1} - x_a^n} $ and $\frac{ {\partial_x{\mathbf A}^n_{\perp,h}(x_a^{n+1})  - \partial_x{\mathbf A}^n_{\perp,h}(x_a^n) } }{x_a^{n+1} - x_a^n} $  are in the form of $\frac{0}{0}$, which could be avoided by  
$$\frac{ {{\mathbf A}^n_{\perp,h}(x_a^{n+1})  - {\mathbf A}^n_{\perp,h}(x_a^n) } }{x_a^{n+1} - x_a^n} \approx \partial_x {\mathbf A}^n_{\perp,h}(\frac{x_a^{n+1} + x_a^n}{2}), \frac{ {\partial_x{\mathbf A}^n_{\perp,h}(x_a^{n+1})  - \partial_x{\mathbf A}^n_{\perp,h}(x_a^n) } }{x_a^{n+1} - x_a^n} \approx \partial_x^2 {\mathbf A}^n_{\perp,h}(\frac{x_a^{n+1} + x_a^n}{2}).$$
\end{remark}
\begin{remark}
Multiplying $\mathbb{G}^{\top}\mathbb{M}_1$ from left with the scheme about ${\mathbf e}_x$, we have 
\begin{equation}
\begin{aligned}
\mathbb{G}^{\top}\mathbb{M}_1{\mathbf e}_x^{n+1} & = \mathbb{G}^{\top}\mathbb{M}_1{\mathbf e}_x^{n}  - \mathbb{G}^{\top}\sum_a \int_{t^n}^{t^{n+1}}{\Lambda}_1(x_a(\tau))^{\top} \mathrm{d}\tau w_a \frac{\mathrm{d}{x_a}(\tau)}{\mathrm{d}\tau},\\
& =  \mathbb{G}^{\top}\mathbb{M}_1{\mathbf e}_x^{n} -  \sum_a \int_{t^n}^{t^{n+1}} \frac{\mathrm{d}}{\mathrm{d}\tau} {\Lambda}_0(x_a(\tau))  w_a \mathrm{d}\tau,\\
& =   \mathbb{G}^{\top}\mathbb{M}_1{\mathbf e}_x^{n} -  {\mathbb{\Lambda}_0({{\mathbf X}^{n+1}})^{\top}} \mathbb{W} \mathbb{1}_{N_p} +  {\mathbb{\Lambda}_0({{\mathbf X}^{n}})^{\top}} \mathbb{W} \mathbb{1}_{N_p},
\end{aligned}
\end{equation}
where $\mathbb{1}_{N_p}$ the vector of size $N_p$ composed of $1$.
Then, the discrete Poisson equation (weak formulation) $\mathbb{G}^\top \mathbb{M}_1 {\mathbf e}_x(t) = - { \mathbb{\Lambda}_0({\mathbf X})^\top} \mathbb{W} \mathbb{1}_{N_p}$ is always satisfied by the numerical solution if it holds initially.
\end{remark}

\noindent{\bf Subsystem II}
The second subsystem about ${\mathbf e}_y,  {\mathbf e}_z,  {\mathbf a}_y,  {\mathbf a}_z$ is 
\begin{equation}
\begin{aligned}
&\dot{\mathbf{X}} = {\mathbf 0},\ \dot{\mathbf P} =  {\mathbf 0},\ \dot{\mathbf S} =  {\mathbf 0},\\
&\dot{\mathbf e}_x = {\mathbf 0},\ \dot{\mathbf e}_y = \mathbb{M}_0^{-1}\left(\sum_{a=1}^{N_p} \frac{\omega_a}{ \sqrt{1+p_a^2 + |{\mathbf A}_\perp(x_a)|^2}}{\Lambda}_0(x_a){A}_{y,h} + \mathbb{G}^{\top}\mathbb{M}_1\mathbb{G}{\mathbf a}_y  \right) - {\textgoth{h}} \, \mathbb{M}_0^{-1}\mathbb{G}^{\top}\mathbb{\Lambda}_{1}({\mathbf X})^{\top}\mathbb{W}{\mathbf{S}}_z,\\
& \gamma(x_a) = \sqrt{1 + p_a^2 + |{\mathbf A}_\perp(x_a)|^2}\\
& \dot{\mathbf e}_z =\mathbb{M}_0^{-1}\left(\sum_{a=1}^{N_p} \frac{\omega_a}{\sqrt{1+p_a^2 + |{\mathbf A}_\perp(x_a)|^2}} {\Lambda}_0(x_a){A}_{z,h} + \mathbb{G}^{\top}\mathbb{M}_1\mathbb{G}{\mathbf a}_z  \right) + {\textgoth{h}} \,\mathbb{M}_0^{-1}\mathbb{G}^{\top}\mathbb{\Lambda}_{1}({\mathbf X})^{\top}\mathbb{W}{\mathbf{S}}_y,\\
&\dot{\mathbf a}_y = -{\mathbf e}_y,\ \dot{\mathbf a}_z = -{\mathbf e}_z.
\end{aligned}
\end{equation}
The discrete gradients about ${\mathbf e}_y,  {\mathbf e}_z,  {\mathbf a}_y,  {\mathbf a}_z$ are
{\small{
$$
\bar{\nabla}_{{\mathbf a}_y} H = \sum_a w_a \frac{({\mathbf A}_{y,h}^{n+1}(x_a^n) + {\mathbf A}_{y,h}^n(x_a^n))\Lambda_1(x_a^n) }{\sqrt{1 + (p_a^n)^2 + |{\mathbf A}^n_{\perp,h}(x_a^n)|^2  } + \sqrt{1 + (p_a^n)^2 + |{\mathbf A}^{n+1}_{\perp,h}(x_a^n)|^2  }  } +  \mathbb{G}^{\top}\mathbb{M}_1\mathbb{G}\frac{ {\mathbf a}_y^n + {\mathbf a}_y^{n +1}}{2} - {\textgoth{h}} \, \mathbb{G}^{\top}\mathbb{\Lambda}_{1}({\mathbf X}^n)^{\top}\mathbb{W}{\mathbf{S}}_z^n,
$$
$$
\bar{\nabla}_{{\mathbf a}_z}H = \sum_a w_a \frac{({\mathbf A}_{z,h}^{n+1}(x_a^n) + {\mathbf A}_{z,h}^n(x_a^n)) \Lambda_1(x_a^n)  }{\sqrt{1 + (p_a^n)^2 + |{\mathbf A}^n_{\perp,h}(x_a^n)|^2  } + \sqrt{1 + (p_a^n)^2 + |{\mathbf A}^{n+1}_{\perp,h}(x_a^n)|^2  }  } +  \mathbb{G}^{\top}\mathbb{M}_1\mathbb{G}\frac{ {\mathbf a}_z^n + {\mathbf a}_z^{n +1}}{2}  + {\textgoth{h}} \, \mathbb{G}^{\top}\mathbb{\Lambda}_{1}({\mathbf X}^n)^{\top}\mathbb{W}{\mathbf{S}}_y^n.
$$
$$
\bar{\nabla}_{{\mathbf e}_y} H = \mathbb{M}_0 \frac{{\mathbf e}_y^{n} + {\mathbf e}_y^{n+1}   }{2}, \quad \bar{\nabla}_{{\mathbf e}_z} H = \mathbb{M}_0 \frac{{\mathbf e}_z^{n} + {\mathbf e}_z^{n+1}   }{2}.
$$
}}
With the above discrete gradient, we have the following scheme,
\begin{equation}
\begin{aligned}
& {\mathbf{X}}^{n+1} =  {\mathbf{X}}^{n},\ {\mathbf P}^{n+1} =  {\mathbf P}^n,\ {\mathbf S}^{n+1} =  {\mathbf S}^n, \ {\mathbf e}_x^{n+1}
 = {\mathbf e}_x^n\\
& \frac{{\mathbf e}_y^{n+1}  - {\mathbf e}_y^{n}}{\Delta t} = \mathbb{M}_0^{-1}\left(\sum_{a=1}^{N_p} \frac{\omega_a({ A}_{y,h}^{n+1}(x_a^n) + { A}_{y,h}^n(x_a^n))\Lambda_1(x_a^n) }{\sqrt{1 + (p_a^n)^2 + |{\mathbf A}^n_{\perp,h}(x_a^n)|^2  } + \sqrt{1 + (p_a^n)^2 + |{\mathbf A}^{n+1}_{\perp,h}(x_a^n)|^2  }  }  + \mathbb{G}^{\top}\mathbb{M}_1\mathbb{G} \frac{{\mathbf a}_y^n +{\mathbf a}_y^{n+1} }{2}  \right) \\
& - {\textgoth{h}} \, \mathbb{M}_0^{-1}\mathbb{G}^{\top}\mathbb{\Lambda}_{1}({\mathbf X}^n)^{\top}\mathbb{W}{\mathbf{S}}_z^n,\\
& \frac{{\mathbf e}_z^{n+1}  - {\mathbf e}_z^{n}}{\Delta t} =\mathbb{M}_0^{-1}\left(\sum_{a=1}^{N_p} \frac{\omega_a({ A}_{z,h}^{n+1}(x_a^n) + { A}_{z,h}^n(x_a^n))\Lambda_1(x_a^n) }{\sqrt{1 + (p_a^n)^2 + |{\mathbf A}^n_{\perp,h}(x_a^n)|^2  } + \sqrt{1 + (p_a^n)^2 + |{\mathbf A}^{n+1}_{\perp,h}(x_a^n)|^2  }  }  + \mathbb{G}^{\top}\mathbb{M}_1\mathbb{G}\frac{{\mathbf a}_z^n +{\mathbf a}_z^{n+1} }{2}  \right) \\
& + {\textgoth{h}} \,\mathbb{M}_0^{-1}\mathbb{G}^{\top}\mathbb{\Lambda}_{1}({\mathbf X}^n)^{\top}\mathbb{W}{\mathbf{S}}_y^n,\\
& \frac{{\mathbf a}_y^{n+1}  - {\mathbf a}_y^{n}}{\Delta t} = -\frac{{\mathbf e}_y^n +{\mathbf e}_y^{n+1} }{2},\  \frac{{\mathbf a}_z^{n+1}  - {\mathbf a}_z^{n}}{\Delta t} = -\frac{{\mathbf e}_z^n +{\mathbf e}_z^{n+1} }{2}.
\end{aligned}
\end{equation}
Then we have 
\begin{equation}
\begin{aligned}
\left( \mathbb{M}_0 + \frac{\Delta t^2}{4} \mathbb{G}^\top \mathbb{M}_1 \mathbb{G}\right) {\mathbf e}_y^{n+1} & = \left( \mathbb{M}_0 - \frac{\Delta t^2}{4} \mathbb{G}^\top \mathbb{M}_1 \mathbb{G}\right) {\mathbf e}_y^{n} + \Delta t \mathbb{G}^{\top}\mathbb{M}_1\mathbb{G} {\mathbf a}_y^n - \Delta t {\textgoth{h}} \, \mathbb{M}_0^{-1}\mathbb{G}^{\top}\mathbb{\Lambda}_{1}({\mathbf X}^n)^{\top}\mathbb{W}{\mathbf{S}}_z \\
& + \Delta t \sum_{a=1}^{N_p} \frac{\omega_a({ A}_{y,h}^{n+1}(x_a^n) + { A}_{y,h}^n(x_a^n))\Lambda_1(x_a^n) }{\sqrt{1 + (p_a^n)^2 + |{\mathbf A}^n_{\perp,h}(x_a^n)|^2  } + \sqrt{1 + (p_a^n)^2 + |{\mathbf A}^{n+1}_{\perp,h}(x_a^n)|^2  }  },\\
\left( \mathbb{M}_0 + \frac{\Delta t^2}{4} \mathbb{G}^\top \mathbb{M}_1 \mathbb{G}\right) {\mathbf e}_z^{n+1} & = \left( \mathbb{M}_0 - \frac{\Delta t^2}{4} \mathbb{G}^\top \mathbb{M}_1 \mathbb{G}\right) {\mathbf e}_z^{n} + \Delta t \mathbb{G}^{\top}\mathbb{M}_1\mathbb{G} {\mathbf a}_z^n + \Delta t {\textgoth{h}} \, \mathbb{M}_0^{-1}\mathbb{G}^{\top}\mathbb{\Lambda}_{1}({\mathbf X}^n)^{\top}\mathbb{W}{\mathbf{S}}_y \\
& + \Delta t \sum_{a=1}^{N_p} \frac{\omega_a({ A}_{z,h}^{n+1}(x_a^n) + { A}_{z,h}^n(x_a^n))\Lambda_1(x_a^n) }{\sqrt{1 + (p_a^n)^2 + |{\mathbf A}^n_{\perp,h}(x_a^n)|^2  } + \sqrt{1 + (p_a^n)^2 + |{\mathbf A}^{n+1}_{\perp,h}(x_a^n)|^2  }  }, 
\end{aligned}
\end{equation}
where on the right side $A_{y,h}^{n+1}, A_{z,h}^{n+1}$ are represented with ${\mathbf e}^n_y, {\mathbf e}^{n+1}_y, {\mathbf e}^n_z, {\mathbf e}^{n+1}_z$ using the equation $  \frac{{\mathbf a}_y^{n+1}  - {\mathbf a}_y^{n}}{\Delta t} = -\frac{{\mathbf e}_y^n +{\mathbf e}_y^{n+1} }{2},\  \frac{{\mathbf a}_z^{n+1}  - {\mathbf a}_z^{n}}{\Delta t} = -\frac{{\mathbf e}_z^n +{\mathbf e}_z^{n+1} }{2}$.
To solve the above scheme about ${\mathbf e}_y^{n+1}, {\mathbf e}_z^{n+1}$, a fixed point iteration method is used combined with a pre-conditioner of $\mathbb{M}_0$. During each iteration, to compute the terms containing ${\mathbf A}^{n+1}_{\perp,h}(x_a^n), 1 \le a \le N_p$ on the right hand side, a loop of all the particles is required.

\noindent{\bf Subsystem III}
The third subsystem is
\begin{equation}
\begin{aligned}
&\dot{\mathbf{X}} = {\mathbf 0},\ \dot{\mathbf P} = {\mathbf 0},\ \dot{\mathbf S} = { \frac{1}{{\textgoth{h}}}}\mathbb{S}\frac{\partial H_{s}}{\partial {\mathbf S}},\\
&\dot{\mathbf e}_x = {\mathbf 0},\ \dot{\mathbf e}_y = {\mathbf 0},\ \dot{\mathbf e}_z = {\mathbf 0},\\
&\dot{\mathbf a}_y = {\mathbf 0},\ \dot{\mathbf a}_z = {\mathbf 0}.
\end{aligned}
\end{equation}
As Hamiltonian depends on ${\mathbf S}$ linearly, discrete gradient for ${\mathbf S}$ is just usual gradient, i.e., $
\bar{\nabla}_{{\mathbf S}}H = {\nabla}_{{\mathbf S}}H $.
For the $a$-th particle, we have
\begin{align}
\label{sdot_hs}
\begin{aligned}
\dot{\mathbf s}_a =
& \left(
\begin{matrix}
    \dot{s}_{a,x}  \\
    \dot{s}_{a,y}  \\
      \dot{s}_{a,z}
      \end{matrix}
\right) = {{
\left(
\begin{matrix}
    0 & Y_a& Z_a \\
    -Y_a  & 0  & 0 \\
      -Z_a & 0 &0
      \end{matrix}
\right)}}
 \left(
\begin{matrix}
    {s}_{a,x}  \\
    {s}_{a,y}  \\
     {s}_{a,z}
      \end{matrix}
\right) =: {{\hat{r}_a}} {\mathbf s}_a,
\end{aligned}
\end{align}
where {{$Y_a = ({\mathbf a}_y^n)^{\top}\mathbb{G}^{\top}{ \Lambda}^1(x_a^n)$, $Z_a = ({\mathbf a}_z^n)^{\top}\mathbb{G}^{\top}{\Lambda}^1(x_a^n)$.}} The Rodrigues' formula gives the following explicit solution for  \eqref{sdot_hs}
\begin{equation}
\label{spin_rodrigue}
{\mathbf s}_a^{n+1} = {{\exp(\Delta t\hat{r}_a){\mathbf s}_a(t^n) = \left(I + \frac{\sin(\Delta t|{\mathbf r}_a|)}{|\mathbf {r}_a|}\hat{r}_a + \frac{1}{2}\left( \frac{\sin(\frac{\Delta t}{2}|{\mathbf r}_a|)}{\frac{|{\mathbf r}_a|}{2}} \right)^2 \hat{r}_a^2\right) {\mathbf s}_a^n}},
\end{equation}
where ${\mathbf r}_a  = (0, Z_a, -Y_a)^\top \in \mathbb{R}^3$, and $I$ is the $3\times 3$ identity matrix.

\subsection{Two dimensional case}
The Poisson matrix~\eqref{eq:dispoisson2d} is split into the following four parts,
\begin{align*}
\begin{aligned}
& \mathbb{J}_1({\mathbf{u}}) \!=\!\!
\left(
\begin{matrix}
    {\mathbf{0}} \!\! & {{\mathbb{W}}}^{-1}  &{\mathbf{0}} & {\mathbf{0}} \!\! & {\mathbf{0}} & {\mathbf{0}}& {\mathbf{0}}  \\
     -{{\mathbb{W}}}^{-1}  \!\! &  {\mathbf{0}}   & {\mathbf{0}} & \mathbb{\Lambda}_1({\mathbf {X}}){\mathbb{M}}_1^{-1} \!\! & {\mathbf{0}} & {\mathbf{0}} & {\mathbf{0}} \\
      {\mathbf{0}}  \!\! & {\mathbf{0}}  &  {\mathbf{0}} & {\mathbf{0}} \!\!& {\mathbf{0}}  & {\mathbf{0}} & {\mathbf{0}} \\
       {\mathbf{0}}  \!\! & -{{\mathbb{M}}}_1^{-1}\mathbb{\Lambda}_1({\mathbf{X}})^{\top}  & {\mathbf{0}}& {\mathbf{0}}  \!\!& {\mathbf{0}}& {\mathbf{0}}& {\mathbf{0}}\\
       {\mathbf{0}}\!\! & {\mathbf{0}}& {\mathbf{0}}& {\mathbf{0}} \!\!& {\mathbf{0}}& {\mathbf{0}} &  {\mathbf{0}}\\
       {\mathbf{0}}\!\! &{\mathbf{0}}&{\mathbf{0}}&{\mathbf{0}}\!\!&{\mathbf{0}}&{\mathbf{0}} & {\mathbf{0}} \\
       {\mathbf{0}}\!\! &{\mathbf{0}}&{\mathbf{0}}&{\mathbf{0}}\!\!&  {\mathbf{0}} & {\mathbf{0}} &{\mathbf{0}}
               \end{matrix}
\right),
\end{aligned}
\begin{aligned}
& \mathbb{J}_2({\mathbf{u}}) \!=\!\!
\left(
\begin{matrix}
    {\mathbf{0}} \!\! & {\mathbf{0}}   &{\mathbf{0}} & {\mathbf{0}} \!\! & {\mathbf{0}} & {\mathbf{0}}& {\mathbf{0}}  \\
      {\mathbf{0}}  \!\! & \mathbb{S}^p  & {\mathbf{0}} & {\mathbf{0}} \!\! & {\mathbf{0}} & {\mathbf{0}} & {\mathbf{0}} \\
      {\mathbf{0}}  \!\! & {\mathbf{0}}  & \frac{1}{\textgoth{h}} {\mathbb{S}} & {\mathbf{0}} \!\!& {\mathbf{0}}  & {\mathbf{0}} & {\mathbf{0}} \\
       {\mathbf{0}}  \!\! & {\mathbf{0}}  & {\mathbf{0}}& {\mathbf{0}}  \!\!& {\mathbf{0}}& {\mathbf{0}}& {\mathbf{0}}\\
       {\mathbf{0}}\!\! & {\mathbf{0}}& {\mathbf{0}}& {\mathbf{0}}  \!\!& {\mathbf{0}}& {\mathbf{0}} &  {\mathbf{0}}\\
       {\mathbf{0}}\!\! &{\mathbf{0}}&{\mathbf{0}}&{\mathbf{0}}\!\!&{\mathbf{0}}&{\mathbf{0}} & {\mathbf{0}} \\
       {\mathbf{0}}\!\! &{\mathbf{0}}&{\mathbf{0}}&{\mathbf{0}}\!\!&  {\mathbf{0}} & {\mathbf{0}} &{\mathbf{0}}
               \end{matrix}
\right),
\end{aligned}\\
\begin{aligned}
& \mathbb{J}_3({\mathbf{u}}) \!=\!\!
\left(
\begin{matrix}
    {\mathbf{0}} \!\! & {\mathbf{0}}  &{\mathbf{0}} & {\mathbf{0}} \!\! & {\mathbf{0}} & {\mathbf{0}}& {\mathbf{0}}  \\
     {\mathbf{0}}  \!\! & {\mathbf{0}}   & {\mathbf{0}} &  {\mathbf{0}}  \!\! & {\mathbf{0}} & {\mathbf{0}} & {\mathbf{0}} \\
      {\mathbf{0}}  \!\! & {\mathbf{0}}  & {\mathbf{0}}  & {\mathbf{0}} \!\!& {\mathbf{0}}  & {\mathbf{0}} & {\mathbf{0}} \\
       {\mathbf{0}}  \!\! & {\mathbf{0}}  & {\mathbf{0}}& {\mathbf{0}}  \!\!& {\mathbf{0}} & {\mathbf{0}}& {\mathbf{0}}\\
       {\mathbf{0}}\!\! & {\mathbf{0}}& {\mathbf{0}}&  {\mathbf{0}}  \!\!& {\mathbf{0}}& {\mathbf{0}} &  {\mathbf{0}}\\
       {\mathbf{0}}\!\! &{\mathbf{0}}&{\mathbf{0}}&{\mathbf{0}}\!\!&{\mathbf{0}}&{\mathbf{0}} &{{\mathbb{M}}}_0^{-1}\\
       {\mathbf{0}}\!\! &{\mathbf{0}}&{\mathbf{0}}&{\mathbf{0}}\!\!&  {\mathbf{0}} &-{{\mathbb{M}}}_0^{-1}&{\mathbf{0}}
               \end{matrix}
\right),
\end{aligned}
\begin{aligned}
& \mathbb{J}_4({\mathbf{u}}) \!=\!\!
\left(
\begin{matrix}
    {\mathbf{0}} \!\! & {\mathbf{0}} &{\mathbf{0}} & {\mathbf{0}} \!\! & {\mathbf{0}} & {\mathbf{0}}& {\mathbf{0}}  \\
     {\mathbf{0}}   \!\! & {\mathbf{0}} & {\mathbf{0}} & {\mathbf{0}} \!\! & {\mathbf{0}} & {\mathbf{0}} & {\mathbf{0}} \\
      {\mathbf{0}}  \!\! & {\mathbf{0}}  & {\mathbf{0}} & {\mathbf{0}} \!\!& {\mathbf{0}}  & {\mathbf{0}} & {\mathbf{0}} \\
       {\mathbf{0}}  \!\! & {\mathbf{0}}  & {\mathbf{0}}& {\mathbf{0}}  \!\!& {\mathbb{M}}_1^{-1}{\mathbb{C}}^\top& {\mathbf{0}}& {\mathbf{0}}\\
       {\mathbf{0}}\!\! & {\mathbf{0}}& {\mathbf{0}}& -{\mathbb{C}} {{\mathbb{M}}}_1^{-1} \!\!& {\mathbf{0}}& {\mathbf{0}} &  {\mathbf{0}}\\
       {\mathbf{0}}\!\! &{\mathbf{0}}&{\mathbf{0}}&{\mathbf{0}}\!\!&{\mathbf{0}}&{\mathbf{0}} & {\mathbf{0}} \\
       {\mathbf{0}}\!\! &{\mathbf{0}}&{\mathbf{0}}&{\mathbf{0}}\!\!&  {\mathbf{0}} & {\mathbf{0}} &{\mathbf{0}}
               \end{matrix}
\right),
\end{aligned}
\end{align*}
which correspond to the following four subsystems.

\noindent{\bf Subsystem I}
The first subsystem about ${\mathbf x}_a, {\mathbf p}_a, {\mathbf e}_{xy}$ is 
\begin{equation}
\begin{aligned}
&\dot{\mathbf x}_a = \frac{{\mathbf p}_a}{\sqrt{1+|{\mathbf p}_a|^2 + |{A}_z({\mathbf x}_a)|^2}},\\
& \dot{\mathbf p}_a = {\mathbf E}_{xy, h}({\mathbf x}_a) - \frac{{{ A}_{z, h}}({\mathbf x}_a)   \nabla_{\mathbf x} { A}_{z, h}({\mathbf x}_a)} {\sqrt{1+|{\mathbf p}_a|^2 + |{A}_z({\mathbf x}_a)|^2}} \\
&  - {\textgoth{h}} s^a_x \left( \partial^2_{x_1x_2} A_z({\mathbf x}_{a}),    \partial^2_{x_2} A_z({\mathbf x}_{a}) \right)^\top + {\textgoth{h}} s^a_y\left(   \partial^2_{x_1} A_z({\mathbf x}_{a}),   \partial^2_{x_1x_2} A_z({\mathbf x}_{a}) \right)^\top - {\textgoth{h}} s^a_z \left(  \partial_{x_1} B_z({\mathbf x}_{a}),    \partial_{x_2} B_z({\mathbf x}_{a}) \right)^\top,\\
&\dot{\mathbf e}_{xy} = -\mathbb{M}_1^{-1} \sum_{a=1}^{N_p} {\Lambda}_1({\mathbf x}_a) w_a  \frac{{\mathbf p}_a}{\sqrt{1+|{\mathbf p}_a|^2 + |{A}_z({\mathbf x}_a)|^2}},\\
&\dot{\mathbf S} =  {\mathbf 0},\  \dot{\mathbf e}_z = {\mathbf 0},\ \dot{\mathbf a}_z = {\mathbf 0}, \ \dot{\mathbf b}_z = {\mathbf 0}.
\end{aligned}
\end{equation}
The discrete gradients of ${\mathbf x}_a, {\mathbf p}_a, {\mathbf e}_{xy}$ are
\begin{equation}
\begin{aligned}
&\bar{\nabla}_{{\mathbf x}_a} H =  w_a\frac{ \left( {{A}^n_{z,h}(x_a^{n+1})  + {A}^n_{z,h}(x_a^n) } \right) \cdot \left( \frac{{{A}^n_{z,h}(x_{1, a}^{n+1}, x_{2, a}^{n+1} )  - { A}^n_{z,h}(x_{1, a}^n, x_{2, a}^{n+1} )}}{(x_{1,a}^{n+1} - x_{1,a}^n)},  \frac{{{A}^n_{z,h}(x_{1, a}^{n}, x_{2, a}^{n+1} )  - { A}^n_{z,h}(x_{1, a}^n, x_{2, a}^{n} )}}{(x_{2,a}^{n+1} - x_{2,a}^n)} \right)^\top    }{\sqrt{ 1 + |{\mathbf p}_a^n|^2 + |{A}^n_{z,h}(x_a^{n+1})|^2 } + \sqrt{ 1 + |{\mathbf p}_a^{n}|^2 + |{A}^n_{z,h}(x_a^{n})|^2 } },\\
& + w_a {\textgoth{h}} s_{a,x}^n \left( \frac{\partial_y A^n_{z,h}(x_{1,a}^{n+1}, x_{2,a}^{n+1} ) -  \partial_y A^n_{z,h}(x_{1,a}^{n}, x_{2,a}^{n+1}) }{x_{1,a}^{n+1} - x_{1,a}^n},     \frac{\partial_y A^n_{z,h}(x_{1,a}^{n}, x_{2,a}^{n+1} ) -  \partial_y A^n_{z,h}(x_{1,a}^{n}, x_{2,a}^{n}) }{x_{2,a}^{n+1} - x_{2,a}^n} \right)^\top,\\
& - w_a {\textgoth{h}} s_{a,y}^n \left( \frac{\partial_x A^n_{z,h}(x_{1,a}^{n+1}, x_{2,a}^{n+1} ) -  \partial_x A^n_{z,h}(x_{1,a}^{n}, x_{2,a}^{n+1}) }{x_{1,a}^{n+1} - x_{1,a}^n},     \frac{\partial_x A^n_{z,h}(x_{1,a}^{n}, x_{2,a}^{n+1} ) -  \partial_x A^n_{z,h}(x_{1,a}^{n}, x_{2,a}^{n}) }{x_{2,a}^{n+1} - x_{2,a}^n} \right)^\top,\\
& + w_a {\textgoth{h}} s_{a,z}^n \left( \frac{B^n_{z,h}(x_{1,a}^{n+1}, x_{2,a}^{n+1} ) -  B^n_{z,h}(x_{1,a}^{n}, x_{2,a}^{n+1}) }{x_{1,a}^{n+1} - x_{1,a}^n},     \frac{B^n_{z,h}(x_{1,a}^{n}, x_{2,a}^{n+1} ) - B^n_{z,h}(x_{1,a}^{n}, x_{2,a}^{n}) }{x_{2,a}^{n+1} - x_{2,a}^n} \right)^\top,\\
&\bar{\nabla}_{{\mathbf p}_a} H = w_a \frac{ {\mathbf p}_a^{n} + {\mathbf p}_a^{n+1}}{\sqrt{ 1 + |{\mathbf p}_a^n|^2 + |{A}^n_{z,h}({\mathbf x}_a^{n+1})|^2 } + \sqrt{ 1 + |{\mathbf p}_a^{n+1}|^2 + |{A}^n_{z,h}({\mathbf x}_a^{n+1})|^2 }}\\
&\bar{\nabla}_{{\mathbf e}_{xy}} H =\mathbb{M}_1  \frac{{\mathbf e}_{xy}^n + {\mathbf e}_{xy}^{n+1} }{2}.
\end{aligned}
\end{equation}
Then we have the following scheme,
\begin{equation}
\begin{aligned}
&\frac{{\mathbf x}_a^{n+1} - {\mathbf x}_a^{n} }{\Delta t} = \frac{ {\mathbf p}_a^{n} + {\mathbf p}_a^{n+1}}{\sqrt{ 1 + |{\mathbf p}_a^n|^2 + |{A}^n_{z,h}({\mathbf x}_a^{n+1})|^2 } + \sqrt{ 1 + |{\mathbf p}_a^{n+1}|^2 + |{A}^n_{z,h}({\mathbf x}_a^{n+1})|^2 }},\\
&\frac{{\mathbf p}_a^{n+1} - {\mathbf p}_a^{n} }{\Delta t} = \frac{1}{\Delta t}\int_{t^n}^{t^{n+1}} {\Lambda}_1^\top({\mathbf x}_a(\tau))\mathrm{d}\tau {\mathbf e}_{xy}  - \frac{1}{w_a}\bar{\nabla}_{{\mathbf x}_a}{H},\\
&\frac{{\mathbf e}_{xy}^{n+1} - {\mathbf e}_{xy}^{n} }{\Delta t}  = -\mathbb{M}_1^{-1} \sum_{a=1}^{N_p}\frac{1}{\Delta t}\int_{t^n}^{t^{n+1}} {\Lambda}_1({\mathbf x}_a(\tau))\mathrm{d}\tau w_a \frac{ {\mathbf p}_a^{n} + {\mathbf p}_a^{n+1}}{\sqrt{ 1 + |{\mathbf p}_a^n|^2 + |{A}^n_{z,h}({\mathbf x}_a^{n+1})|^2 } + \sqrt{ 1 + |{\mathbf p}_a^{n+1}|^2 + |{A}^n_{z,h}({\mathbf x}_a^{n+1})|^2 }},\\
&{\mathbf S}^{n+1} =  {\mathbf S}^{n},\  {\mathbf e}_z^{n+1} =  {\mathbf e}_z^{n},\, {\mathbf a}_z^{n+1} =  {\mathbf a}_z^{n} ,\ {\mathbf b}_z^{n+1} =  {\mathbf b}_z^{n},
\end{aligned}
\end{equation}
where the time-continuous trajectory is defined as 
$$
{\mathbf x}_a(\tau) = {\mathbf x}_a^n + (\tau - t^n) \frac{{\mathbf x}_a^{n+1} - {\mathbf x}_a^n}{\Delta t}, \quad \tau \in [t^n, t^{n+1}], \quad 1 \le a \le N_p.
$$
Similar to remark~\ref{re:dispoisson}, we can also prove discrete Poisson equation is satisfied by the numerical solution.

\noindent{\bf Subsystem II}
The second subsystem is 
$$
\dot{\mathbf p}_a = (  \frac{p_{a,y} B^n_{z,h}({\mathbf x}_a)}{\sqrt{1 + |{\mathbf p}_a|^2 + A_{z,h}^2({\mathbf x}_a)} },  - \frac{p_{a,x} B_z^n({\mathbf x}_a)}{\sqrt{1 + |{\mathbf p}_a|^2 + A_{z,h}^2({\mathbf x}_a)} }  )^\top, \quad \dot{\mathbf s}_a = {\mathbf s}_a \times {\mathbf B}^n_h({\mathbf x}_a), \quad 1 \le a \le N_p.
$$
As $|{\mathbf p}_a|^2$ is conserved by this subsystem, 
the velocity and spin variables can be solved exactly using Rodrigues' formula as~\eqref{spin_rodrigue}, and naturally energy is conserved.

\noindent{\bf Subsystem III}
The third subsystem is 
\begin{equation}
\begin{aligned}
\dot{{\mathbf a}_z} &= - {\mathbf e}_z,\\
 \dot{ {\mathbf e}_z} &= \mathbb{M}_0^{-1}\nabla_{{\mathbf a}_z}H =  \mathbb{M}_0^{-1}\left(\mathbb{G}_{*}^\top \mathbb{M}_{1,*} \mathbb{G}_{*} {\mathbf a}_z + {\textgoth{h}} \mathbb{G}_{*}^\top \mathbb{\Lambda}_{1,*}({\mathbf X})^\top \mathbb{W} {\mathbf S}_{xy}  \right)\\
 & +  \mathbb{M}_0^{-1} \sum_a w_a \frac{{ A}_{z,h}({\mathbf x}_a)\Lambda_0({\mathbf x}_a) }{\sqrt{1+|{\mathbf p}_a|^2 + |{A_z}({\mathbf x}_a)|^2}}.
\end{aligned}
\end{equation}
With the discrete gradient about ${\mathbf e}_z$ and ${\mathbf a}_z$, 
\begin{equation}
\begin{aligned}
\bar{\nabla}_{{\mathbf e}_z}H &= \frac{ {\mathbf e}_z^n +  {\mathbf e}_z^{n+1} }{2},\\
\bar{\nabla}_{{\mathbf a}_z}H & = \mathbb{G}_{*}^\top \mathbb{M}_{1,*} \mathbb{G}_{*}  \frac{ {\mathbf a}_z^n +  {\mathbf a}_z^{n+1} }{2} +  {\textgoth{h}} \mathbb{G}_{*}^\top \mathbb{\Lambda}_{1,*}({\mathbf X}^n)^\top \mathbb{W} {\mathbf S}_{xy}^n \\
& + \sum_a w_a \frac{({ A}_{z,h}^n({\mathbf x}_a^n)+{ A}_{z,h}^{n+1}({\mathbf x}_a^n))\Lambda_0({\mathbf x}_a) }{\sqrt{1+|{\mathbf p}_a^n|^2 + |{A_z^n}({\mathbf x}_a^n)|^2}+\sqrt{1+|{\mathbf p}_a^n|^2 + |{A_z^{n+1}}({\mathbf x}_a^n)|^2}},
\end{aligned}
\end{equation}
we have the following scheme,
\begin{equation}
\begin{aligned}
 \frac{ {\mathbf a}_z^{n+1} -  {\mathbf a}_z^{n} }{\Delta t} &= - \frac{ {\mathbf e}_z^n +  {\mathbf e}_z^{n+1} }{2},\\
 \frac{ {\mathbf e}_z^{n+1} -  {\mathbf e}_z^{n} }{\Delta t} &= \mathbb{M}_0^{-1}\nabla_{{\mathbf a}_z}H =  \mathbb{M}_0^{-1}\left(  \mathbb{G}_{*}^\top \mathbb{M}_{1,*} \mathbb{G}_{*}  \frac{ {\mathbf a}_z^n +  {\mathbf a}_z^{n+1} }{2} +  {\textgoth{h}} \mathbb{G}_{*}^\top \mathbb{\Lambda}_{1,*}({\mathbf X}^n)^\top \mathbb{W} {\mathbf S}_{xy}^n\right)\\
 & +  \mathbb{M}_0^{-1}  \sum_a w_a \frac{({ A}_{z,h}^n({\mathbf x}_a^n)+{ A}_{z,h}^{n+1}({\mathbf x}_a^n))\Lambda_0({\mathbf x}_a^n) }{\sqrt{1+|{\mathbf p}_a^n|^2 + |{A_z^n}({\mathbf x}_a^n)|^2}+\sqrt{1+|{\mathbf p}_a^n|^2 + |{A_z^{n+1}}({\mathbf x}_a^n)|^2}}.
\end{aligned}
\end{equation}
After substituting the above first equation into second one, we get 
\begin{equation}
\begin{aligned}
\left( \mathbb{M}_0 + \frac{\Delta t^2}{4}  \mathbb{G}_{*}^\top \mathbb{M}_{1,*} \mathbb{G}_{*}  \right) {\mathbf e}_z^{n+1} &= \left( \mathbb{M}_0- \frac{\Delta t^2}{4}  \mathbb{G}_{*}^\top \mathbb{M}_{1,*} \mathbb{G}_{*}  \right) {\mathbf e}_z^{n} + \Delta t  \mathbb{G}_{*}^\top \mathbb{M}_{1,*} \mathbb{G}_{*}  {\mathbf a}_z^n +   {\textgoth{h}} \Delta t \mathbb{G}_{*}^\top \mathbb{\Lambda}_{1,*}({\mathbf X}^n)^\top \mathbb{W} {\mathbf S}_{xy}\\
& + \Delta t \sum_a w_a \frac{({ A}_{z,h}^n({\mathbf x}_a^n)+{ A}_{z,h}^{n+1}({\mathbf x}_a^n))\Lambda_0({\mathbf x}_a^n) }{\sqrt{1+|{\mathbf p}_a^n|^2 + |{A_z^n}({\mathbf x}_a^n)|^2}+\sqrt{1+|{\mathbf p}_a|^2 + |{A_z^{n+1}}({\mathbf x}_a^n)|^2}},
\end{aligned}
\end{equation}
where on the right side $A_{z,h}^{n+1}$ is represented with ${\mathbf e}^n_z, {\mathbf e}^{n+1}_z$ using the equation $ \frac{ {\mathbf a}_z^{n+1} -  {\mathbf a}_z^{n} }{\Delta t} = - \frac{ {\mathbf e}_z^n +  {\mathbf e}_z^{n+1} }{2}$.
The above equation about ${\mathbf e}_z^{n+1}$ can be solved with the fixed point iteration method combined with a pre-conditioner of $\mathbb{M}_0$. In each iteration, a loop of all the particles is required to compute the terms related with ${A_z^{n+1}}({\mathbf x}_a^n), 1 \le a \le N_p$.

\noindent{\bf Subsystem IV}
The fourth subsystem is 
\begin{equation}
\begin{aligned}
\dot{{\mathbf e}_{xy}} &= \mathbb{M}_1^{-1}\mathbb{C}^\top \nabla_{{\mathbf b}_z}H = \mathbb{M}_1^{-1}\mathbb{C}^\top\left( \mathbb{M}_2 {\mathbf b}_z +   {\textgoth{h}}  \mathbb{\Lambda}_2({\mathbf X})^\top \mathbb{W} {\mathbf S}_z \right),\\
 \dot{ {\mathbf b}_z} &= - \mathbb{C} \mathbb{M}_1^{-1} \nabla_{{\mathbf e}_{xy}}H = - \mathbb{C} {\mathbf e}_{xy}.
\end{aligned}
\end{equation}
With the discrete gradient about ${\mathbf e}_{xy}$ and ${\mathbf b}_z$, 
\begin{equation}
\begin{aligned}
\bar{\nabla}_{{\mathbf e}_{xy}}H &= \frac{ {\mathbf e}_{xy}^n +  {\mathbf e}_{xy}^{n+1} }{2},\\
\bar{\nabla}_{{\mathbf b}_z}H & = \mathbb{M}_{2}  \frac{ {\mathbf b}_{z}^n +  {\mathbf b}_{z}^{n+1} }{2} +   {\textgoth{h}}  \mathbb{\Lambda}_2({\mathbf X}^n)^\top \mathbb{W} {\mathbf S}_z^n,
\end{aligned}
\end{equation}
we have the following scheme,
\begin{equation}
\begin{aligned}
 \frac{ {\mathbf e}_{xy}^{n+1} -  {\mathbf e}_{xy}^{n} }{\Delta t} &= \mathbb{M}_1^{-1}\mathbb{C}^\top\left(  \mathbb{M}_{2}  \frac{ {\mathbf b}_{z}^n +  {\mathbf b}_{z}^{n+1} }{2} +   {\textgoth{h}}  \mathbb{\Lambda}_2({\mathbf X}^n)^\top \mathbb{W} {\mathbf S}_z^n \right),\\
 \frac{ {\mathbf b}_z^{n+1} -  {\mathbf b}_z^{n} }{\Delta t} &= - \mathbb{C} \frac{ {\mathbf e}_{xy}^n +  {\mathbf e}_{xy}^{n+1} }{2},
\end{aligned}
\end{equation}
from which we get 
\begin{equation}
\begin{aligned}
\left(\mathbb{M}_1 + \frac{\Delta t^2}{4}\mathbb{C}^\top \mathbb{M}_2 \mathbb{C} \right) {\mathbf e}_{xy}^{n+1}  = \left(\mathbb{M}_1 - \frac{\Delta t^2}{4}\mathbb{C}^\top \mathbb{M}_2 \mathbb{C} \right) {\mathbf e}_{xy}^{n}  + \Delta t \mathbb{C}^\top \mathbb{M}_2 {\mathbf b}_z^n + \Delta t  {\textgoth{h}} \mathbb{C}^\top \mathbb{\Lambda}_2({\mathbf X}^n)^\top \mathbb{W} {\mathbf S}_z^n. 
\end{aligned}
\end{equation}
A fixed point iteration method is used to solve above equation with a pre-conditioner for $\mathbb{M}_1$.

\section{Numerical experiments}
In this section, two one dimensional numerical experiments are done for two cases: without spin effect and with spin effects. In both cases, energy conservation property is validated, also we found that discrete Poisson equation is satisfied by the numerical solution indeed. Moreover, in the former case, numerical growth rates of Fourier modes are compared with the analytical ones. In both cases, iteration tolerance is set as $10^{-13}$, and B-spline degrees in $V_0$ and $V_1$ are 3 and 2, respectively.
\subsection{Without spin effect}
\begin{figure}[htbp]
\center{
\subfigure[]{\includegraphics[scale=0.36]{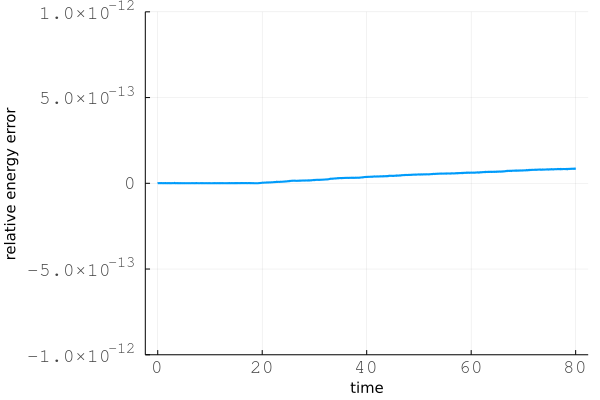}}
\subfigure[]{\includegraphics[scale=0.36]{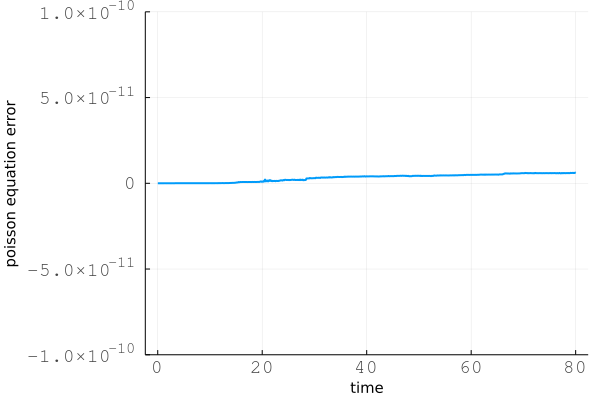}}
}
\caption{{\bf Without spin} Time evolutions of relative energy error and poisson equation error.}
\label{fig:without1}
\end{figure}

\begin{figure}[htbp]
\center{
\subfigure[]{\includegraphics[scale=0.36]{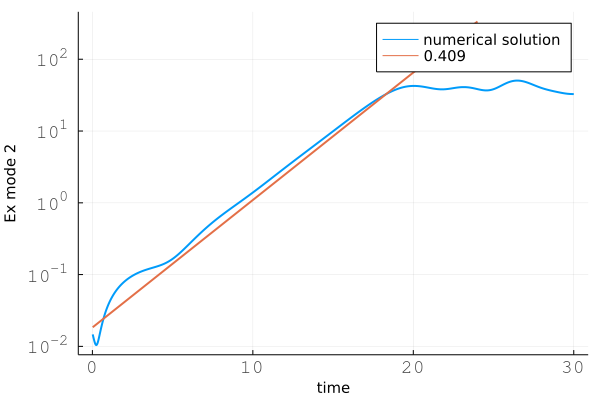}}
\subfigure[]{\includegraphics[scale=0.36]{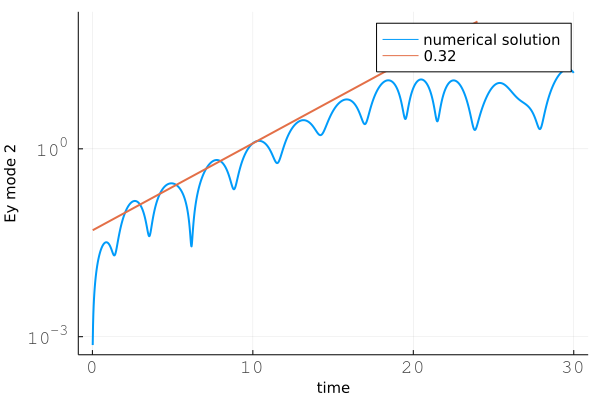}}
}
\caption{{\bf Without spin} Time evolutions of the amplitude of the second Fourier mode of $E_x$ and the amplitude of the second Fourier mode of $E_y$.}
\label{fig:without2}
\end{figure}

In this numerical test, which is called parametric instability~\cite{D1}, ${\textgoth{h}} = 0$, i.e., spin effects is not included, we take the initial distribution function as a homogeneous Maxwellian expressed as
$$
f_0(x, p) = \frac{1}{\sqrt{2\pi T}} \exp(-\frac{p^2}{2T}), \quad T = 3/511.
$$
The initial conditions of fields are
$$
E_{x0}(x) = 0, E_{y0} = E_0 \cos(kx), E_{z0}(x) = E_0 \sin(kx), A_{y0}(x) = -E_0 \sin(kx), A_{z0} = E_0 \cos(kx),
$$
where $E_0 = \sqrt{3}, k = \frac{1}{\sqrt{2}}$. Simulation space domain is $[0, \frac{2\pi}{k}]$, time step size is $\Delta t = 0.02$, final simulation time is 80, cell number in space is 128, particle number is $4 \times 10^3$, and Lie--Trotter splitting is used. 
The time evolutions of relative energy error and Poisson equation error are plotted in Fig.~\ref{fig:without1}. We can see that the error is at the level of iteration tolerance, and has no obvious growth with time. In Fig.~\ref{fig:without2}, we compare the numerical growth rates of the second Fourier mode of $E_x$ and $E_y$ with analytical rates (red lines)~\cite{D1}, which fit in well and validate the code.

\subsection{With spin effects}
\begin{figure}[htbp]
\center{
\subfigure[]{\includegraphics[scale=0.36]{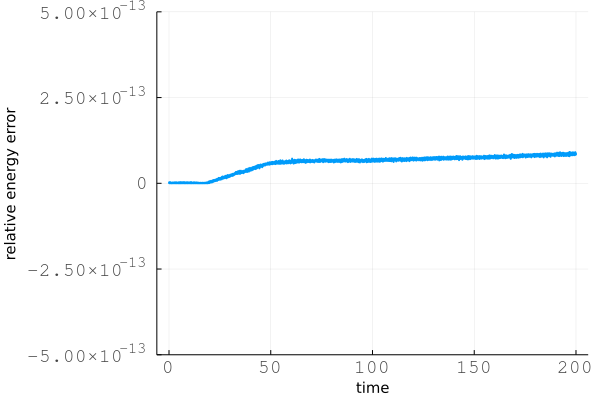}}
\subfigure[]{\includegraphics[scale=0.36]{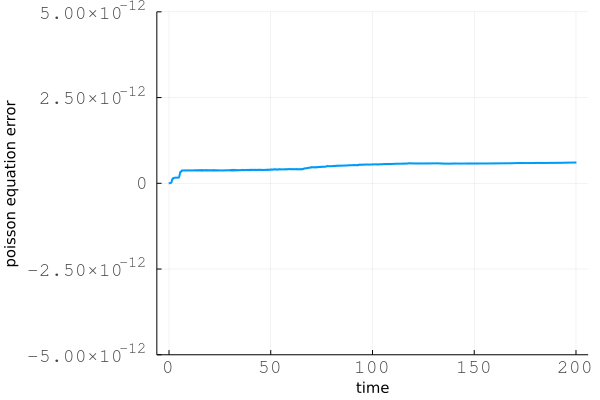}}
}
\caption{{\bf With spin} Time evolutions of relative energy error and poisson equation error.}
\label{fig:with1}
\end{figure}

\begin{figure}[htbp]
\center{
\subfigure[]{\includegraphics[scale=0.36]{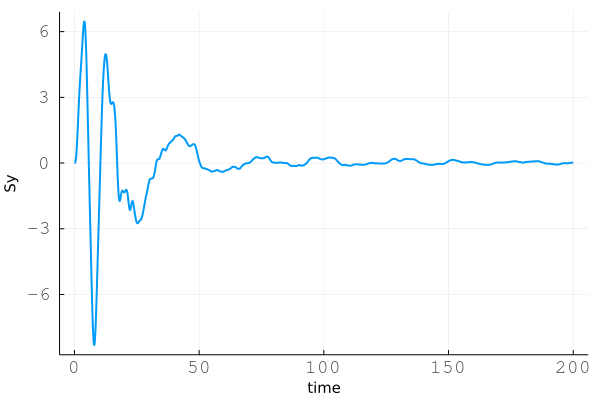}}
\subfigure[]{\includegraphics[scale=0.36]{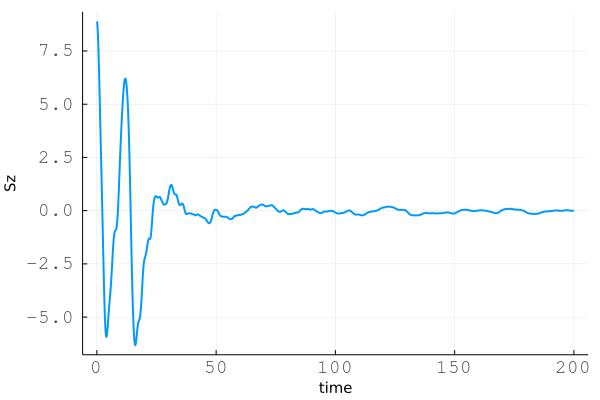}}
}
\caption{{\bf With spin} Time evolutions of the $S_y = \int s_y f \mathrm{d}{\mathbf s}\mathrm{d}{p}\mathrm{d}x$ and $S_z = \int s_z f \mathrm{d}{\mathbf s}\mathrm{d}{p}\mathrm{d}x$.}
\label{fig:with2}
\end{figure}

In this numerical test, we include spin effects by setting ${\textgoth{h}} = 0.1$. We take the same initial condition as the case without spin effects but with a different initial distribution function, i.e.,
$$
f_0(x, p, {\mathbf s}) = \frac{1}{\sqrt{2\pi T}} \exp(-\frac{p^2}{2T})\mathbb{1}_{\{ 1\}}({\mathbf s}), \quad \mathbb{1}_{\{ 1\}}({\mathbf s}) = \left\{
\begin{aligned}
  &   1, \quad {\mathbf s} = (0, 0, 1)^\top\\
 &0, \quad \text{else} \\
  \end{aligned}\right., \quad T = 3/511,
$$
$$
E_{x0}(x) = 0, E_{y0} = E_0 \cos(kx), E_{z0}(x) = E_0 \sin(kx), A_{y0}(x) = -E_0 \sin(kx), A_{z0} = E_0 \cos(kx),
$$
where $E_0 = \sqrt{3}, k = \frac{1}{\sqrt{2}}$. Simulation space domain is $[0, \frac{2\pi}{k}]$, time step size is $\Delta t  = 0.02$, final simulation time is 200, cell number in space is 128, particle number is $10^4$, and Lie--Trotter splitting is used. 
From Fig.~\ref{fig:with1}, we can see that the energy error and poisson equation error are quite small and have no obvious growth with time. In Fig.~\ref{fig:with2}, time evolution of spin momentum at $y$ and $z$ directions are plotted, we find that the momentums oscillate with time and decay to zeros finally, which are similar to the results of non-relativistic case in~\cite{LPFEEC}.

\section{Conclusion}
In this work, discrete gradient method is used to construct energy conserving particle-in-cell schemes for one and two dimensional relativistic Vlasov--Maxwell equations with spin effects. The space discretization of fields is done in the framework of finite element exterior calculus. Numerical experiments are done to validate our numerical schemes, especially the conservation properties. Three dimension case is not detailed in this work, as the relativistic factor $\sqrt{1+|{\mathbf p}|^2}$ does not depends on particle position, and thus is easier to apply the discrete gradient method. There are several future works to be envisaged,  such as parallelization could be done to accelerate the code, non-periodic boundary condition can also be considered as~\cite{boundary} to conduct more practical simulations.

\section{Appendix}
\subsection{Three dimensional spin Vlasov--Maxwell equations}
\label{sec:threed}
\begin{equation}
\label{eq:Vlasovn_norm}
\begin{aligned}
\frac{\partial f}{\partial t}&+\frac{\bf p}{\gamma}\cdot \nabla f+[\left({\bf E}+\frac{\bf p}{\gamma}\times{\bf B}\right) + \textgoth{h} \; \nabla({\bf s} \cdot {\bf B})]\cdot\frac{\partial f}{\partial{\bf p}} + ({\bf s}\times {\bf B}) \cdot \frac{\partial f}{\partial {\bf s}}=0,\\
\frac{\partial{\bf E}}{\partial t}&= \nabla\times{\bf B}- \int_{\mathbb{R}^6} \frac{\bf p}{\gamma} f \mathrm{d}{\mathbf p}\mathrm{d}{\mathbf s} +  \textgoth{h} \; \nabla \times \int_{\mathbb{R}^6}  {\mathbf s} f \mathrm{d}{\mathbf p}\mathrm{d}{\mathbf s}, \\
\frac{\partial{\bf B}}{\partial t} &= - \nabla\times{\bf E},\\
\nabla\cdot{\bf E}&={\int_{\mathbb{R}^6} f\mathrm{d}{\mathbf p}\mathrm{d}{\mathbf s} -1}, \\
\nabla\cdot{\bf B}&=0,
\end{aligned}
\end{equation}
where ${\mathbf x} = (x_1, x_2, x_3)^\top \in \mathbb{R}^3, {\mathbf p} = (p_x, p_y, p_z)^\top \in \mathbb{R}^3, {\mathbf s} = (s_x, s_y, s_z)^\top \in \mathbb{R}^3$.
Hamiltonian of the above system is 
$$
\mathcal{H} = \int (\sqrt{1+|{\mathbf p}|^2} - 1 ) f \mathrm{d}{\mathbf x}  \mathrm{d}{\mathbf p}  \mathrm{d}{\mathbf s}  + \frac{1}{2}\int |{\mathbf E}|^2  \mathrm{d}{\mathbf x}  +    \frac{1}{2}\int |{\mathbf B}|^2  \mathrm{d}{\mathbf x} + \textgoth{h} \int {\mathbf s} \cdot {\mathbf B} f \mathrm{d}{\mathbf x}  \mathrm{d}{\mathbf p}  \mathrm{d}{\mathbf s}. 
$$

\subsection{Discrete functional derivatives of 2D reduced model}
\begin{equation}
\begin{aligned}
&\frac{\delta \mathcal{F}}{\delta E_z} = (\Lambda^0)^\top \mathbb{M}_0^{-1} \nabla_{{\mathbf e}_z}F,  \quad \frac{\delta \mathcal{F}}{\delta A_z} = (\Lambda^0)^\top \mathbb{M}_0^{-1} \nabla_{{\mathbf a}_z}F, \\
& \frac{\delta \mathcal{F}}{\delta E_{xy}} = (\Lambda^1)^\top \mathbb{M}_1^{-1} \nabla_{{\mathbf e}_{xy}}F,\quad \frac{\delta \mathcal{F}}{\delta B_z} = (\Lambda^2)^\top \mathbb{M}_2^{-1} \nabla_{{\mathbf b}_z}F, \\
& \frac{\partial}{\partial {\mathbf x}}\frac{\delta \mathcal{F}}{\delta f}|({\mathbf x}_a,  {\mathbf p}_a, {\mathbf s}_a) = \frac{1}{w_a} \nabla_{{\mathbf x}_a}F, \quad \frac{\partial}{\partial {\mathbf p}}\frac{\delta \mathcal{F}}{\delta f}|({\mathbf x}_a,  {\mathbf p}_a, {\mathbf s}_a) = \frac{1}{w_a} \nabla_{{\mathbf p}_a}F, \frac{\partial}{\partial {\mathbf s}}\frac{\delta \mathcal{F}}{\delta f}|({\mathbf x}_a,  {\mathbf p}_a, {\mathbf s}_a) = \frac{1}{w_a} \nabla_{{\mathbf s}_a}F.
\end{aligned}
\end{equation}

\label{sec:dis2d}

\end{document}